\documentclass[a4paper,11pt]{amsart} 
\usepackage{amsmath,amsxtra,amssymb,latexsym, amscd,amsthm}
\usepackage[mathscr]{eucal}
\usepackage{mathrsfs}
\usepackage{bbm}
\usepackage{enumerate}
\usepackage{pict2e}
\usepackage{tikz}
\usepackage{graphicx}
\usepackage[a4paper]{geometry}
\usepackage{color}
\usepackage{hyperref}
\numberwithin{equation}{section}
\theoremstyle{plain}
\newtheorem{thm}{Theorem}[section]
\newtheorem{prp}[thm]{Proposition}
\newtheorem{cor}[thm]{Corollary}

\newtheorem*{euc*}{Euclidean division}
\newtheorem*{fek*}{Fekete's Lemma}
\newtheorem*{kin*}{Kingman's Subadditive Ergodic Theorem}
\newtheorem*{fur*}{Furstenberg-Kesten Theorem}
\theoremstyle{definition}

\newtheorem{pro}[thm]{Problem}
\newtheorem{rem}[thm]{Remark}

\newtheorem*{rem*}{Remark}

\newcommand{\dd}{\mathrm{d}}

 \newcommand{\To}{\longrightarrow}
\renewcommand{\Im}{\operatorname{Im}}
\renewcommand{\Re}{\operatorname{Re}}

\makeatletter
\newcommand*{\ov}[1]{%
  $\m@th\overline{\mbox{#1}}$%
}
\newcommand*{\ovA}[1]{%
  $\m@th\overline{\mbox{#1}\raisebox{3mm}{}}$%
}
\newcommand*{\ovB}[1]{%
  $\m@th\overline{\mbox{#1\rule{0pt}{3mm}}}$%
}
\newcommand*{\ovC}[1]{%
  $\m@th\overline{\mbox{#1\strut}}$%
}
\newcommand*{\ovD}[1]{%
  $\m@th\overline{\mbox{#1\vphantom{\"A}}}$%
}
\newcommand*{\ovE}[1]{%
  $\m@th\overline{\raisebox{0pt}[1.2\height]{#1}}$%
}
\newcommand*{\ovF}[1]{%
  $\m@th\overline{\raisebox{0pt}[\dimexpr\height+0.3mm\relax]{#1}}$%
}
\newcommand*{\ovG}[1]{%
  $\m@th\overline{\raisebox{0pt}[\dimexpr\height+1mm\relax]{#1\vphantom{A}}}$%
}
\makeatother
\newcommand{\N}{\mathbb{N}}
\newcommand{\Z}{\mathbb{Z}}

\newcommand{\R}{\mathbb{R}}
\newcommand{\C}{\mathbb{C}}

\newcommand{\T}{\mathbb{T}}

\newcommand{\eps}{\epsilon}
\newcommand{\cW}{\mathcal{W}}

\newcommand{\cH}{\mathcal{H}}
\newcommand{\cG}{\mathcal{G}}

\newcommand{\cA}{\mathcal{A}}
\newcommand{\cB}{\mathcal{B}}
\newcommand{\cT}{\mathcal{T}}
\newcommand{\cM}{\mathcal{M}}
\newcommand{\cN}{\mathcal{N}}
\newcommand{\cS}{\mathcal{S}}

\newcommand{\tilg}{{\tilde{g}}}

\DeclareMathOperator{\expect}{\mathbb{E}}

\DeclareMathOperator{\prob}{\mathbb{P}}

\DeclareMathOperator{\supp}{supp}

\DeclareMathOperator{\varnbi}{Var_{nb}^I}

\DeclareMathOperator{\var}{Var\,}

\DeclareMathOperator{\vari}{Var^I}

\DeclareSymbolFont{extraup}{U}{zavm}{m}{n}
\DeclareMathSymbol{\varheart}{\mathalpha}{extraup}{86}
\DeclareMathSymbol{\vardiamond}{\mathalpha}{extraup}{87}

\title{Recent results of quantum ergodicity on graphs and further investigation}
\author{Nalini Anantharaman and Mostafa \textsc{Sabri}}
\address{Universit\'e de Strasbourg, CNRS, IRMA UMR 7501, F-67000 Strasbourg, France.}
\email{anantharaman@math.unistra.fr}
\address{Department of Mathematics, Faculty of Science, Cairo University, Cairo 12613, Egypt.}
\address{Universit\'e Paris Sud XI, UMR 8628 du CNRS, Laboratoire de Math\'ematique, B\^at. 307, 91405 Orsay Cedex, France.}
\email{mostafa.sabri@math.u-psud.fr}
\usepackage{calc}
\usepackage{yhmath}
\usepackage{graphicx}

\makeatletter
\newlength{\temp@wc@width}
\newlength{\temp@wc@height}
\newcommand{\widecheck}[1]{%
  \setlength{\temp@wc@width}{\widthof{$#1$}}%
  \setlength{\temp@wc@height}{\heightof{$#1$}}%
  #1\hspace{-\temp@wc@width}%
  \raisebox{\temp@wc@height+2pt}[\heightof{$\widehat{#1}$}]%
     {\rotatebox[origin=c]{180}{\vbox to 0pt{\hbox{$\widehat{\hphantom{#1}}$}}}}%
}
\makeatother

\begin{document}

\begin{abstract}
We outline some recent proofs of quantum ergodicity on large graphs and give new applications in the context of irregular graphs. We also discuss some remaining questions.
\end{abstract}

\maketitle

\section{Introduction}          \label{sec:introd}

Our aim in this paper is to present a comprehensive outline of the main ideas behind some recent proofs of quantum ergodicity on large graphs \cite{A,AS2} or \cite{BLL}. We hope this will be a helpful guide for an audience with different backgrounds~: spectral geometry, quantum chaos, graph theory and Anderson localization. We also include new applications to our main results and present some open problems. 

In its original context, the quantum ergodicity theorem, also dubbed Shnirelman theorem, asserts that the eigenfunctions of the Laplace-Beltrami operator on a compact manifold become uniformly distributed in the large eigenvalue limit, provided the geodesic flow is ergodic \cite{Shni,Zel, CdV}.

In the framework of combinatorial graphs, quantum ergodicity is interpreted in the large spatial limit; that is, we consider finite graphs whose size goes to infinity. 
Quantum ergodicity is a result of \emph{spatial delocalization} for the eigenfunctions of the adjacency matrix (and more generally, discrete Schr\"odinger operators) on a large graph. The first result of this kind, proven by Anantharaman and Le Masson \cite{ALM}, asserts that if a sequence $(G_N)$ of finite graphs ``converges'' to the $(q+1)$-regular tree $\T_q$ and if the graphs are expanders, then most eigenfunctions of $\cA_{G_N}$ become uniformly distributed on $G_N$ when $N$ gets large. This somehow complements the phenomenon of \emph{spectral delocalization}, meaning that the $\ell^2$-spectrum of the adjacency matrix on the tree $\T_q$ is purely absolutely continuous.

We should clarify what we mean by the ``convergence'' of $(G_N)$ to $\T_q$. An adequate notion is the \emph{Benjamini-Schramm} or  \emph{local weak} convergence \cite{BS}. This is a convergence in a distributional sense. One may define a distance on the set of rooted graphs by asserting that two rooted graphs are close to each other if there are large isomorphic balls around their respective roots. The local weak convergence of $(G_N)$ is not the convergence in this metric, but rather the weak convergence of the associated probability measures $(U_{G_N})$, defined by choosing a root uniformly at random in $G_N$. Saying that $(G_N)$ converges to $\T_q$ means more precisely that in the space of probability measures on rooted graphs, $U_{G_N}$ converges weakly to the Dirac mass $\delta_{[\T_q,o]}$, where $o\in \T_q$ is an arbitrary root. Such convergence holds iff for any $r\in \N$, $\frac{\# \{x\in V_N:B_{G_N}(x,r) \cong B_{\T_q}(o,r)\}}{|V_N|} \To 1$ as $N\to \infty$. Here we assumed the $G_N$ are connected, as we do throughout the paper, and denoted $B_G(x,r)$ the $r$-ball centered at $x$ in $G$.

As we mentioned above, quantum ergodicity can be regarded as a spatial delocalization phenomenon, counterpart of the fact that the $\ell^2$-spectrum $\sigma(\cA_{\T_q})$ is reduced to the absolutely continuous spectrum $\sigma_{ac}(\cA_{\T_q})$. Suppose then more generally, that we have a (possibly random) Schr\"odinger operator $\cH$ on some (possibly random) tree $\cT$ which has purely absolutely continuous spectrum. Suppose that we are given a (deterministic) sequence of discrete Schr\"odinger operators on graphs $(G_N)$, converging to $(\cT, \cH)$ in the local weak sense. Does quantum ergodicity hold for the eigenfunctions of $H_{G_N}$ on $G_N$ as $N$ gets large~? This is the question we studied in \cite{AS2} and to which we provided a positive answer, if in addition the graphs are expanders. This generalized the results of Anantharaman--Le Masson \cite{ALM} to non-regular graphs, and to discrete Schr\"odinger operators more general than the adjacency matrix. As an important application \cite{AS3}, we proved the quantum ergodicity phenomenon for the Anderson model on $\T_q$, in the regime of weak disorder. We recall that spectral delocalization for the Anderson model on $\T_q$ at weak disorder was first proven by Klein \cite{Klein} (with some precursor ideas by Kunz and Souillard \cite{KuSo}). Our paper \cite{AS3} thus completes the picture by proving spatial delocalization. This is one of very few delocalization results for the Anderson model, the regimes of localization (pure point spectrum, exponential decay of eigenfunctions) being rather well understood today. Shedding some light on delocalization for models on \emph{euclidean} graphs such as $\Z^d$, $d\ge 3$, remains to this day a major open problem.

Our quantum ergodicity theorems are stated for {\em{deterministic}} sequences of graphs, in particular, they hold for explicitly defined families of expanders. One could also ask what happens for ``typical'' (i.e. random) sequences. Very strong delocalization results have recently been proved by Erd\"os, Knowles, Yau and Yin \cite{EKYY} for the eigenfunctions of the adjacency matrix of Erd\"os--R\'enyi graphs, and by Bauerschmidt, Huang and Yau \cite{BHY} for random regular graphs. Note that these results hold for \emph{almost all} graphs, but do not say what happens for given deterministic examples. A more detailed comparison between random and deterministic results is provided in our paper \cite{AS2}. 

The paper is organized as follows. In Section~\ref{sec:adj}, we give an essentially complete proof of quantum ergodicity for the adjacency matrix on regular graphs. In Section~\ref{sec:schro} we discuss the general case of Schr\"odinger operators on graphs of bounded degree. New applications of our abstract result appear in Section~\ref{sec:cone}, where we consider spacial delocalization in trees of finite cone type. Throughout, we leave some open questions.

\section{Quantum ergodicity for the adjacency matrix of regular graphs}\label{sec:adj}

\subsection{The result}

Let $(G_N)$ be a sequence of $(q+1)$-regular graphs, $q\ge 2$, say $G_N=(V_N,E_N)$ with $|V_N|=N$. We assume $(G_N)$ converges to $\T_q$ in the local weak sense. This holds iff $G_N$ has few short cycles, more precisely, $\frac{\# \{x:\rho_{G_N}(x)<r\}}{N}\to 0$ for any $r$, where $\rho_{G_N}(x)$ is the injectivity radius at $x\in V_N$ -- the largest $\rho$ such that $B_{G_N}(x,\rho)$ is a tree. We refer to this assumption as \textbf{(BST)}.

The aim is to prove that the eigenfunctions of $\cA_{G_N}$ become delocalized when $N$ gets large. We shall assume that $(G_N)$ is a sequence of expanders, that is, there is $\beta>0$ such that the spectrum $\sigma(\frac{\cA_{G_N}}{q+1})$ is contained in $[-1+\beta,1-\beta]\cup \{1\}$ for all $N$. We call this assumption \textbf{(EXP)}.

Both assumptions are generic, in the sense that random sequences of $(q+1)$-regular graphs converge to $\T_q$ almost surely (as follows from Corollary 2.19 in Bollob\'as's book \cite{Bol}), and they are expanders with probability tending to one; the latter fact is often attributed to Pinsker, however, to find a full proof we refer to \S 6.2 in \cite{ACF}. Friedman reinforced this, proving that random regular graphs are almost Ramanujan with probability tending to one \cite{Fri}.
There are also explicit sequences of $(q+1)$-regular graphs satisfying both assumptions, for instance the Ramanujan graphs constructed by Lubotzky-Phillips-Sarnak \cite{LPS}, or Cayley graphs of $SL(2, \Z/p\Z)$ studied by Helfgott and Bourgain-Gamburd \cite{BG,Hel}. It is important to note that, for Cayley graphs (or graphs whose automorphism group acts transitively on the vertices), our result does not convey much information if $\psi$ is an eigenfunction associated to a simple eigenvalue~: in this case, $|\psi(x)|$ is going to be independent of the vertex $x$. Fortunately for us, in the Bourgain--Gamburd case, the multiplicity of eigenvalues is very high. Even more fortunately, there are many regular graphs for which the automorphism group  is trivial~: in fact, this is the case for typical $(q+1)$-regular graphs \cite{Bol2}.

We now state the theorem. Assume both \textbf{(BST)} and \textbf{(EXP)} hold for $(G_N)$.


\begin{thm}\label{t:thm1}
Let $(\psi_j^{(N)})_{j=1}^N$ be an orthonormal basis of eigenfunctions of $\cA_{G_N}$ for $\ell^2(V_N)$, with corresponding eigenvalues $(\lambda_j^{(N)})_{j=1}^N$. Fix $R\in\N$ and let ${\mathbf K}_N$ be a $V_N\times V_N$ matrix satisfying ${\mathbf K}_N(x,y)=0$ if $d(x,y)>R$ and $\sup_N\sup_{x,y\in V_N}|{\mathbf K}_N(x,y)|\le 1$. Then
\begin{equation}\label{thm:adj}
\lim_{N\to \infty} \frac{1}{N} \sum_{j=1}^N \left| \langle \psi_j^{(N)},{\mathbf K}_N\psi_j^{(N)}\rangle - \langle {\mathbf K}_N\rangle_{\lambda_j^{(N)}}\right|^2 =0\,,
\end{equation}
\emph{where $\langle {\mathbf K}\rangle_{\lambda} = \frac{1}{N} \sum_{x,y\in V_N} {\mathbf K}(x,y)\Phi(\lambda,d(x,y))$ and}
$\Phi(\lambda, \cdot)$ is the spherical function
\begin{equation}\label{e:spheri}
\Phi(\lambda,r) = q^{-r/2} \left(\frac{2}{q+1}P_r\Big(\frac{\lambda}{2\sqrt{q}}\Big) + \frac{q-1}{q+1} Q_r\Big(\frac{\lambda}{2\sqrt{q}}\Big)\right) .
\end{equation}

Here $P_r(\cos\theta) = \cos(r\theta)$ and $Q_r(\cos \theta)=\frac{\sin(r+1)\theta}{\sin \theta}$ are the Chebyshev polynomials.

\end{thm}

This theorem implies that for $N$ large enough, most $\langle \psi_j,K\psi_j\rangle$ approach $\langle K\rangle_{\lambda_j}$~: for any $\varepsilon >0$,
\begin{equation}\label{e:Markov}
\lim_{N\to \infty} \frac{1}{N} \#\left\{j\in [1, N]: \left| \langle \psi_j^{(N)},K_N\psi_j^{(N)}\rangle - \langle K_N\rangle_{\lambda_j^{(N)}}\right|>\varepsilon\right\}=0\,.
\end{equation}

 In the special case where $R=0$, $K_N=a_N$ is a function on $V_N$ and $\langle \psi_j^{(N)},K_N\psi_j^{(N)}\rangle=\sum_{x\in V_N} a_N(x)| \psi_j^{(N)}(x)|^2$. Also, $\Phi(\lambda, 0)=1$, so $\langle a_N\rangle_{\lambda} $ is simply the uniform average of $a_N$, $\langle a_N\rangle_{\lambda}=\frac{1}{N}\sum_{x\in V_N} a_N(x)$. Taking for instance $a_N = \chi_{\Lambda_N}$, the characteristic function of a set $\Lambda_N\subset V_N$ of size $\alpha N$, $0<\alpha<1$, this implies that for most $\psi_j$, we have $\|\chi_{\Lambda_N}\psi_j\|^2 \approx \alpha$. In particular, if we consider {\em{any}} set containing half the vertices of $V_N$, we will find half the mass of $\|\psi_j\|^2$ in it, \emph{for most $j$}.

We like to interpret this theorem as an ``equidistribution'' result for eigenfunctions, but one should be very cautious with the meaning of this. Note that we are measuring the distance between the uniform measure on $V_N$ and the probability measure $\sum_{x=1}^N| \psi_j^{(N)}(x)|^2\delta_x$
 in a very weak sense, by comparing the average of only one function $a_N$, so \emph{our result is very far from saying that $| \psi_j^{(N)}(x)|^2$ is uniformly close to $\frac1N$}.
 The set of indices $j$ such that $\langle \psi_j^{(N)},a_N\psi_j^{(N)}\rangle$ is close to $\frac{1}{N}\sum_{x\in V_N} a_N(x)$ depends on the test function $a_N$.  Since the speed of convergence we can reach in \eqref{thm:adj} is typically of order $\frac{1}{\log N}$, we can only improve \eqref{e:Markov} by allowing a logarithmic number of observables or a logarithmic error $\varepsilon$, and thus we stay far from comparing $| \psi_j^{(N)}(x)|^2$ and $\frac1N$ at microscopic scale.
 
More generally, the theorem implies that if $d_{G_N}(x,y)=r$, then $\overline{\psi_j(x)}\psi_j(y) \approx \frac{\Phi(\lambda, r)}{N}$ in the same weak sense.

As we mentioned in the introduction, strong delocalization results have been proved in \cite{BHY} in the context of \emph{random} graphs. In particular, \cite{BHY} contains bounds of the form $\|\psi_j\|_{\infty} \lesssim \frac{1}{\sqrt{N}}$ modulo logarithmic corrections, for all $\psi_j$ in the bulk of the spectrum; as well as a result of Quantum Unique Ergodicity, saying that for any given observable $a$, $\sum_{x\in V_N} a(x)| \psi_j^{(N)}(x)|^2$ is close to $\frac{1}{N}\sum_{x\in V_N} a(x)$ for all $\psi_j$, for \emph{most} random regular graphs. This is pretty much the ideal case of uniform distribution of $\psi_j$, where $\psi_j \approx (\frac{1}{\sqrt{N}},\dots,\frac{1}{\sqrt{N}})$. Does one also have the ideal behavior of eigenfunction correlation predicted by \eqref{thm:adj} in this random set-up ?
This extension seems far from obvious, as the method of \cite{BHY} assumes that the test functions are probabilistically independent from the random graph $(G_N)$, and thus cannot be generalized to test observables ${\mathbf K}_N$ that depend on the distance between points in $(G_N)$.

\begin{pro}
For a typical random regular graph $(G_N)$,
prove that for all eigenvectors in the bulk of the spectrum of $\mathcal{A}_{G_N}$, we have $\overline{\psi_j^{(N)}(x)}\psi_j^{(N)}(y) \approx \frac{\Phi(\lambda,r)}{N}$ for all $x\in V_N$ and $y$ such that $d_{G_N}(x,y)=r$.
\end{pro}

\subsection{Spaces of operators \label{s:op}}

Before we present a proof of \eqref{thm:adj}, we introduce some notation.


Fix $G=G_N$, $G=(V,E)$. Let $B_k$ be the set of non-backtracking paths $(x_0,\dots,x_k)$ of length $k$ in $G$. Denote $B=B_1$. We let $\mathscr{H}_k = \C^{B_k}$, the space of complex-valued functions on $B_k$, and $\mathscr{H}_{\le k} = \mathop \oplus_{j=0}^k \mathscr{H}_j$.

For elements of $B_k$, we use the short-hand notation $(x_0;x_k)$ instead of $(x_0,\dots,x_k)$, and for $K\in\mathscr{H}_k$ we write $K(x_0;x_k)$ for $K(x_0,\dots,x_k)$.

Any $K\in \mathscr{H}_{k}$ defines an operator $K_G$ on $\ell^2(V)$ by
\begin{equation}\label{eq:kg}
\langle \varphi, K_G \psi\rangle_{\ell^2(V)} = \sum_{(x_0;x_k)\in B_k} \overline{\varphi(x_0)}K(x_0;x_k)\psi(x_k)
\end{equation}
for $\varphi,\psi \in \ell^2(V)$. In particular, the matrix elements $K_G(x,y) = \langle \delta_x,K_G\delta_y\rangle$ are given by $K_G(x,y) =\sum_{(x_0;x_k)} K(x_0;x_k)$, where the sum runs over all $(x_0;x_k)$ in $B_k$ with fixed endpoints $x_0=x$ and $x_k=y$. We have $K_G(x,y)=0$ unless $x$ and $y$ can be joined by a non-backtracking path of length $ k$. If $\rho_G(x)\ge k$, where $\rho_G(x)$ is the injectivity radius at $x$, observe that the sum defining $K_G(x,y)$ contains at most one nonzero term. 

If $K=(K_j)\in \mathscr{H}_{\le k}$, we extend \eqref{eq:kg} linearly by $\langle \varphi,K_G\psi\rangle = \sum_{j=0}^k \langle \varphi,(K_j)_G\psi\rangle$.

Given $K=(K_j)\in \mathscr{H}_{\le k}$, we consider its normalized $\ell^2$-norm
\begin{equation}\label{hknorm}
\|K\|_{\mathscr{H}}^2 = \sum_{j=0}^k \|K_j\|_{\mathscr{H}_j}^2\qquad \text{where} \quad \|K_j\|_{\mathscr{H}_j}^2 = \frac{1}{N}\sum_{(x_0;x_j)\in B_j} |K(x_0;x_j)|^2 \,.
\end{equation}
We define the Hilbert space $\mathscr{H}$ as the completion of  $\mathop \oplus_{j=0}^{+\infty} \mathscr{H}_j$ for that norm.

We also consider the normalized Hilbert-Schmidt norms of operators,
\[
\|K_G\|_{HSN}^2 = \frac{1}{N}\sum_{x,y\in V_N} |K_G(x,y)|^2 \,.
\]
These two norms are not the same, but they coincide if $K\in \mathscr{H}_{\le k}$ and if the injectivity radius is everywhere greater than $k$.
More generally, we see that for any $K\in \mathscr{H}_{\le k}$,
\begin{equation}\label{e:2norms}
\|K_G\|_{HSN}^2 \le \|K\|_{\mathscr{H}}^2 + c_{k,q} \frac{\# \{x\in V :\rho_G(x) <k\}}{N} \cdot \|K\|_{\infty}^2\,,
\end{equation}
where $\|K\|_{\infty} = \sup_{x,y}|K(x;y)|$ and $c_{k,q} = |B_{G_N}(x,r)|^2 =[1+(q+1) \sum_{j=1}^kq^{j-1}]^2$. Under assumption \textbf{(BST)}, the second term is $o_k(1)_{N\To\infty}$, by which we mean that it depends on the parameter $k$ and goes to $0$ as $N\to \infty$.

In analogy to \eqref{eq:kg}, any $K\in \mathscr{H}_k$ defines an operator $K_B$ on $\ell^2(B)$ by
\begin{equation}\label{e:kb}
\langle f,K_Bg\rangle_{\ell^2(B)} = \sum_{(x_0;x_k) \in B_k} \overline{f(x_0,x_1)}K(x_0;x_k)g(x_{k-1},x_k)
\end{equation}
for $f,g\in \ell^2(B)$. This extends to $K\in \mathscr{H}_{\le k}$. We finally define
\begin{equation}\label{e:mean}
\mathscr{H}_k^o = \{K\in \mathscr{H}_k : \langle K\rangle = 0\} \qquad \text{where} \quad \langle K\rangle = \frac{1}{N}\sum_{(x_0;x_k)\in B_k}K(x_0;x_k)\,.
\end{equation}

\subsection{The quantum variance}\label{sec:qv}
Recall that the spherical function $r\mapsto\Phi(\lambda,r)$ satisfies $\Phi(\lambda,0)=1$, $\Phi(\lambda,1) = \frac{\lambda}{q+1}$ and the linear recursive formula
\begin{equation}\label{eq:ind}
\Phi(\lambda, r+1) = \frac{1}{q} \left(\lambda \Phi(\lambda, r) - \Phi(\lambda, r-1)\right)
\end{equation}
for $r\ge 1$. This is equivalent to the fact that for any $x\in \T_q$, the radial function $y\mapsto \Phi(\lambda, d_{\T_q}(x, y))$ is an eigenfunction of $\cA_{\T_q}$ for the eigenvalue $\lambda$.


Define the constant function $S_k = \frac{1}{(q+1)q^{k-1}} \in \mathscr{H}_k$ for $k\ge 1$, $S_0=1\in\mathscr{H}_0$. Then $[(S_k)_G\psi](x)  = \frac{1}{(q+1)q^{k-1}}\sum_{(x_0;x_k),x_0=x} \psi(x_k)$ for $k\ge 1$ and $[(S_0)_G\psi](x)=\psi(x)$. Note that we have the identity, very specific to regular graphs~:
\begin{equation}\label{e:spher0}
(S_k)_G = \Phi(\cA_G, k)\,,
\end{equation}
where $\Phi(\cA_G, k)$ is defined by replacing $\lambda$ by $\cA_G$ in the polynomial $\lambda\mapsto\Phi(\lambda, k)$. This is easily seen by induction on $k$. In fact, this holds for $k=0,1$ using $\Phi(\lambda,0)=1$ and $\Phi(\lambda,1)=\frac{\lambda}{q+1}$. Next, we have $[(S_{k+1})_G\psi](x)=\frac{1}{(q+1)q^k}\sum_{(x_0;x_{k+1}),x_0=x}\psi(x_{k+1}) = \frac{1}{(q+1)q^k}\sum_{(x_0;x_k),x_0=x}[(\cA_G\psi)(x_k)-\psi(x_{k-1})] = \frac{1}{q}[(S_k\cA_G\psi)(x)-(S_{k-1}\psi)(x)]$. If the identity holds for $j\le k+1$, then using, \eqref{eq:ind}, this becomes $[\Phi(\cA_G,k+1)\psi](x)$, proving \eqref{e:spher0}.

Back to \eqref{thm:adj}, we write
\begin{equation}\label{e:K1}
\langle \psi_j, {\mathbf{K}}\psi_j \rangle - \langle {\mathbf{K}}\rangle_{\lambda_j} = \sum_{k=0}^R \sum_{x_0\in V_N}\sum_{d_{G_N}(x_0,x_k)=k} {\mathbf{K}}(x_0,x_k)\left[\overline{\psi_j(x_0)}\psi_j(x_k) - \frac{\Phi(\lambda_j, k)}{N}\right]\,
\end{equation}
for ${\mathbf{K}}$ an operator with bounded range, as in the assumptions of Theorem \ref{t:thm1}.

On the other hand, we can define $K_k\in \mathscr{H}_k$ by $K_k(x_0;x_k):={\mathbf{K}}(x_0,x_k)$, and $K=(K_k)\in \mathscr{H}_{\leq R}$. We get
\begin{equation}\label{e:K2}
\langle \psi_j, K_G\psi_j \rangle - \sum_{k=0}^R\langle K_k\rangle\Phi(\lambda_j,k) = \sum_{k=0}^R \sum_{(x_0;x_k)\in B_k} K(x_0;x_k)\left[\overline{\psi_j(x_0)}\psi_j(x_k) - \frac{\Phi(\lambda_j, k)}{N}\right]\,
\end{equation}
where $\langle K_k\rangle$ is defined in \eqref{e:mean}.

Note that \eqref{e:K1} and \eqref{e:K2} are almost the same, but the second expression may contain more terms, as there may be several paths $(x_0;x_k)$ with endpoints $x_0,x_k$ if $\rho_G(x_0)<R$. 
But by assumption \textbf{(BST)}, the vertices with $\rho_G(x_0)<R$ are relatively few. Thus if we define the \emph{quantum variance}
\[
\var(K) = \frac{1}{N} \sum_{j=1}^N \left|\langle \psi_j, K_G \psi_j\rangle\right|^2
\]
for $K=(K_k)\in \mathscr{H}_{\leq k}$, then using \eqref{e:spher0}, we get
\[
\frac{1}{N}\sum_{j=1}^N \left| \langle \psi_j,{\mathbf{K}}\psi_j\rangle - \langle {\mathbf{K}}\rangle_{\lambda_j}\right|^2 = \var\Big(K-\sum_{k=0}^R \langle K_k\rangle S_k\Big)  + o_R(1)_{N\To\infty} \,.
\]
By the Cauchy-Schwarz inequality 
\[
\var\Big(K-\sum_{k=0}^R \langle K_k\rangle S_k\Big) \le R\,\sum_{k=0}^R\var(K_k-\langle K_k\rangle S_k) \,.
\]
Finally, note that if $K_k\in\mathscr{H}_k$, then $K_k-\langle K_k\rangle S_k\in \mathscr{H}_k^o$. From this we conclude that in order to prove \eqref{thm:adj} it suffices to show that $\var(K) \to 0$ as $N\to \infty$, for any $K\in \mathscr{H}_k^o$.

\subsection{The proof}\label{sec:proof1}
To show that $\var(K)\to 0$ for any $K\in \mathscr{H}_k^o$, we roughly proceed as follows. Using the eigenfunction equation, we show that the quantum variance possesses some  invariance property, of the form $\var(K) \approx \var(\frac{1}{n}\sum_{r=1}^n \mathcal{C}_r K)$ for some operators $\mathcal{C}_r$. Clearly, we have $\var(K)\le \|K_G\|_{HSN}^2$. So using \eqref{e:2norms}, we get $\var(K) \lesssim \|\frac{1}{n}\sum_{r=1}^n \mathcal{C}_r K\|_{\mathscr{H}}^2$. The goal is to find terms $\mathcal{C}_rK$ that are (almost) orthogonal as $r$ varies, so that this reduces to $\frac{1}{n^2}\sum_{r=1}^n \|\mathcal{C}_rK\|_{\mathscr{H}}^2$. Finally, if we have $\|\mathcal{C}_rK\|_{\mathscr{H}}^2\approx \|K\|_{\mathscr{H}}^2$, the upshot will be that $\var(K) \lesssim \frac{1}{n}\|K\|_{\mathscr{H}}^2$ for any $n$, and the theorem will be proven by taking $n$ arbitrarily large.

To find the operators $\mathcal{C}_r$, note that we have the na\"{i}ve relation
\[
\langle \psi_j K_G\psi_j\rangle = \frac{1}{\lambda_j} \langle \psi_j,K_G \mathcal{A}_G\psi_j\rangle = \frac{1}{n}\sum_{r=1}^n \frac{1}{\lambda_j^r} \langle \psi_j,K_G \mathcal{A}_G^r\psi_j\rangle
\]
at least for $\lambda_j\neq 0$. So one could consider $\mathcal{C}_r = \frac{1}{\lambda_j^r}\cA_G^r$ as a first idea. However, since the powers of $\cA_G$ involve backtracking trajectories, the orthogonality we are seeking for the terms $\mathcal{C}_rK$ is not achieved. An important message of \cite{A} is to consider powers of the non-backtracking adjacency matrix $\mathcal{B}$ instead, which is defined on $\ell^2(B)$ by
\begin{equation}\label{e:nonba}
(\cB f)(x_0,x_1) = \sum_{x_2\in \cN_{x_1}\setminus\{x_0\}} f(x_1,x_2)\,.
\end{equation}
Here $\cN_x$ is the set of nearest neighbours of $x$.

For this to work, the idea is to consider some modified quantum variance in which the invariance law involves only non-backtracking trajectories. This will be implemented in two ways; a very simple one in this section, and a more general one in Section~\ref{sec:schro}.

Given $k\ge 1$, define $\nabla:\mathscr{H}_{k-1}\to\mathscr{H}_k$ by
\[
(\nabla K)(x_0;x_k) = K(x_1;x_k) - K(x_0;x_{k-1})\,.
\]
Its adjoint $\nabla^{\ast}:\mathscr{H}_{k+1}\to\mathscr{H}_k$ is given by
\[
(\nabla^{\ast}K)(x_0;x_k) = \sum_{x_{-1}\in\cN_{x_0}\setminus\{x_1\}} K(x_{-1};x_k) - \sum_{x_{k+1}\in\cN_{x_k}\setminus\{x_{k-1}\}} K(x_0;x_{k+1})
\]
for $k\ge 1$, and for $k=0$, $(\nabla^{\ast}K)(x_0) = \sum_{x_{-1}\sim x_0} K(x_{-1},x_0) - \sum_{x_1\sim x_0} K(x_0,x_1)$.

\medskip

We now show that it suffices to control $\var(\nabla^{\ast}K)$ for any $K\in\mathscr{H}_k^o$. As we shall see, this slight modification of the quantum variance will allow us, in the case of regular graphs, to obtain an invariance law which may be expressed in terms of the non-backtracking random walk (for non-regular graphs, the required modification is less simple).

First, note that for $(q+1)$-regular graphs, $q\ge 2$, the kernel of $\nabla$ is reduced to constants, and thus the range of $\nabla^{\ast}:\mathscr{H}_{k+1}\to\mathscr{H}_k$ coincides with $\mathscr{H}_k^o$. Hence, any $K\in\mathscr{H}_k^o$ may be written in the form $\nabla^{\ast}K'$. The norm of $K'$ in $\mathscr{H}$ is controlled by the spectral gap $\beta$, but in \eqref{e:2norms} we are also using the $\|\cdot\|_{\infty}$-norm, which is not controlled by the spectral gap. For this reason it is better to work with some explicit approximation of $K'$.

Introduce the transfer operator $\cS :\mathscr{H}_k\to\mathscr{H}_k$ by
\[
(\cS K)(x_0;x_k) = \frac{1}{q} \sum_{x_{-1}\in\cN_{x_0}\setminus\{x_1\}} K(x_{-1};x_{k-1})
\]
for $k\ge 1$, and $\cS = \frac{\cA_G}{q+1}$ for $k=0$. Next, define for $T\in\N^{\ast}$,
\[
\cS_TK = \frac{1}{T}\sum_{r=0}^{T-1}(T-r)\cS^rK \quad \text{and} \quad \widetilde{\cS}_TK = \frac{1}{T}\sum_{r=1}^T \cS^rK \,.
\]
Then
\[
K=(I-\cS)(\cS_TK) + \widetilde{\cS}_T K \,.
\]
On the other hand, if $i_k :\mathscr{H}_k\to\mathscr{H}_{k+1}$ is the operator defined by $(i_kK)(x_0;x_{k+1}) = K(x_0;x_k)$, then we see that $\nabla^{\ast}i_kK = q(\cS-I) K $ for $k\ge 1$, and $\nabla^{\ast} i_0K = (q+1)(\cS-I) K $ for $k=0$. Thus (say for $k\geq 1$) we have written
\begin{equation}\label{e:uneequation}
K= -\frac{1}q\nabla^{\ast} i_k\cS_TK+ \widetilde{\cS}_T K \,.
\end{equation}
The operator $-\frac{1}q i_k\cS_T$ is precisely the approximation of $\nabla^{\ast -1}$ we are looking for. In fact, if $K\in \mathscr{H}_k^o$, the remainder term $\widetilde{\cS}_T K$ can be bounded by
\[
\|\widetilde{S}_TK\|_{\mathscr{H}}^2 \le \frac{1}{T^2}\Big(\sum_{r=1}^T \|\cS^r K\|_{\mathscr{H}}\Big)^2 \le \frac{C^2_{k,\beta}}{T^2}\|K\|_{\mathscr{H}}^2
\]
with a constant $C_{k,\beta}$ depending only on $k$ and on the spectral gap $\beta$. This is clear for $k=0$ using \textbf{(EXP)}, since $\|\frac{\cA_G}{q+1}\|_{\mathscr{H}_0^o\to\mathscr{H}_0^o}\le 1-\beta$. It was shown in \cite{ALM} that \textbf{(EXP)} implies an analogous control on $\|\cS^r\|_{\mathscr{H}_k^o\to\mathscr{H}_k^o}$ for all $k$ in case of regular graphs, and this was extended to irregular graphs in \cite{A17}. Thus, $\widetilde{\cS}_T K$ may be made arbitrarily small by choosing $T$ large.

Back to \eqref{e:uneequation}, we get for $K\in \mathscr{H}_k^o$,
\begin{equation}\label{e:red}
\var(K) \le 2q^{-2}\var(\nabla^{\ast}i_k\cS_T K) + 2\var(\widetilde{S}_T K) \,
\end{equation}
and $\var(\widetilde{S}_T K)\lesssim  \|\widetilde{S}_TK\|_{\mathscr{H}}^2  \le \frac{C^2_{k,\beta}}{T^2}\|K\|_{\mathscr{H}}^2$. 
Fixing $T$ large, we are reduced to showing that $\var(\nabla^{\ast}K_T)$ goes to $0$ as $N$ goes to $+\infty$, for $K_T = i_k\cS_T K$. Thus, it suffices to prove we have $\lim_{N\To +\infty}\var(\nabla^{\ast}K) =0$  
for any $K\in \mathscr{H}_k$ to obtain \eqref{thm:adj}. Note that in contrast to the aforementioned $K'=\nabla^{*-1}K$, we have an explicit bound $\|K_T\|_{\infty} \le T\,\|K\|_{\infty}$, so the error term arising from \eqref{e:2norms} is controlled as long as we take $N\to\infty$ before $T\to\infty$.

\medskip

Let $K \in \mathscr{H}_k$, $k\ge 1$. Since $\psi_j$ is an eigenfunction of $\cA_G$, it is clear that for any $j$,
\begin{equation}\label{e:commutator}
\langle \psi_j, [\cA_G, K_G]\psi_j\rangle  =  0 \, ,
\end{equation}
where $[\cdot, \cdot]$ is the commutator of operators. On the other hand, an elementary algebraic calculation reveals that
\[
\langle \cA_G \psi_j, K_G \psi_j\rangle - \langle \psi_j, K_G \cA_G\psi_j\rangle = \langle \psi_j, (\nabla K)_G \psi_j\rangle + \langle \psi_j, (\nabla^{\ast}K)_G\psi_j\rangle\,.
\]
Define $\mathcal{M}^{\ast}:\mathscr{H}_k\to \mathscr{H}_{k+2}$ by $(\mathcal{M}^{\ast}K)(x_0;x_{k+2}) = q^{-1}K(x_1;x_{k+1})$. A specificity of regular graphs is that $-\nabla^{\ast}\cM^{\ast}K = \nabla K$. So \eqref{e:commutator} may be rewritten as
\begin{equation}\label{e:smuggle}
\langle \psi_j, (\nabla^{\ast}\mathcal{M}^{\ast}K)_G\psi_j\rangle  = \langle \psi_j, (\nabla^{\ast}K)_G\psi_j\rangle \, .
\end{equation}

Letting 
\begin{equation}\label{e:Sigma}
\Sigma^n = \frac{1}{n} \sum_{r=1}^n \mathcal{M}^{\ast\, r}\,,
\end{equation}
we get by iteration of \eqref{e:smuggle}
\begin{equation}\label{e:invariantrel}
\var(\nabla^{\ast} \Sigma^n K)  = \var(\nabla^{\ast} K) \, .
\end{equation}
Now note that $\nabla^{\ast} \Sigma^n K = \frac{1}{n} \sum_{r=1}^n \nabla^{\ast} \mathcal{M}^{\ast\, r}K \in \mathop \oplus_{r=1}^n \mathscr{H}_{k+2r-1}$. Using \eqref{e:2norms},
\[
\var(\nabla^{\ast} \Sigma^n K)  \le \|\nabla^{\ast} \Sigma^n K\|_{\mathscr{H}}^2 + c_{q,k,n} \frac{\# \{ x \in V: \rho_G(x) \le k+2n-1 \}}{N} \|\nabla^{\ast} \Sigma^n K\|^2_{\infty} \, .
\]
Now
\[
(\nabla^{\ast} \Sigma^n K)_{k+2r-1}(x_0;x_{k+2r-1}) = \frac{1}{n} (\nabla^{\ast} \mathcal{M}^{\ast\, r} K)(x_0;x_{k+2r-1}) = \frac{-1}{nq^{r-1}} (\nabla K)(x_{r-1};x_{k+r}) \, ,
\]
hence, $\|\nabla^{\ast} \Sigma^n K\|_{\infty} \leq \frac{1}{n} \|\nabla K\|_{\infty} \le \frac{2}{n}\|K\|_{\infty}$. In the Hilbert space $\mathscr{H}$ the various terms $\nabla^{\ast} \mathcal{M}^{\ast\, r}K$ making up
$\nabla^{\ast} \Sigma^n K$ are orthogonal, and thus (by Pythagoras)
\begin{equation}\label{eq:ortho}
\|\nabla^{\ast} \Sigma^n K\|_{\mathscr{H}}^2 = \sum_{r=1}^n \|\nabla^{\ast} \Sigma^n K\|_{\mathscr{H}_{k+2r-1}}^2 = \sum_{r=1}^n\frac{1}{n^2} \|\nabla K\|_{\mathscr{H}_{k+1}}^2 \le \frac{4q}{n} \|K\|_{\mathscr{H}_k}^2\,.
\end{equation}
Recalling \eqref{e:invariantrel}, we have shown that
\[
\var(\nabla^{\ast} K)  \le \frac{4q}{n} \|K\|_{\mathscr{H}_k}^2+ \tilde{c}_{q,k,n} \frac{\# \{ x \in V: \rho_G(x) \le k+2n-1 \}}{N} \|K\|^2_{\infty} \, .
\]
Taking $N\to\infty$ followed by $n\to\infty$, we get $\lim_{N\To +\infty}\var(\nabla^{\ast}K) =0$, completing the proof of \eqref{thm:adj}.

\begin{rem}\label{r:1}
Define $\tau_{\pm} :\ell^2(V)\to \ell^2(B)$ by $(\tau_-\varphi)(x_0,x_1)=\varphi(x_0)$ and $(\tau_+\varphi)(x_0,x_1)=\varphi(x_1)$. We note that, for any $\varphi,\psi\in \ell^2(V)$, $K\in\mathscr{H}_k$, $k\ge 1$,
\begin{equation}\label{e:nabla*}
\langle \varphi, (\nabla^{\ast}K)_G \psi \rangle = \langle \tau_+ \varphi, K_B \tau_+ \psi \rangle - \langle \tau_- \varphi, K_B \tau_- \psi \rangle
\end{equation}
where $K_B : \ell^2(B)\To \ell^2(B)$ was defined in \eqref{e:kb}. Moreover,
\begin{equation}\label{e:m*}
\langle f, (\mathcal{M}^{\ast}K)_B g \rangle = q^{-1} \langle   f, \cB K_B\cB\, g\rangle
\end{equation}
for any $f,g \in \ell^2(B)$. So for any integer $r$, $\langle f, (\mathcal{M}^{\ast \,r}K)_B g \rangle = q^{-r} \langle   f, \cB^r K_B\cB^r\, g\rangle$. Combining \eqref{e:nabla*} with \eqref{e:m*}, we see that the iterations $\nabla^{\ast}\mathcal{M}^{\ast\,r}$ amount to iterations of the non-backtracking random walk $\cB$.

\end{rem}
 
\begin{rem}\label{r:2}
Assuming $\cA_G \psi_j= \lambda_j\psi_j$, one easily checks that
\begin{alignat}{4}\label{e:ultra}
& \cB\,\tau_- \psi_j = q \tau_+ \psi_j \, , & \qquad &  & & \cB\,\tau_+\psi_j = \lambda_j \tau_+ \psi_j - \tau_- \psi_j \, , & \\
& \cB^{\ast} \tau_- \psi_j = \lambda_j \tau_- \psi_j - \tau_+ \psi_j & \qquad &\text{and} &\qquad &\cB^{\ast} \tau_+ \psi_j = q \tau_- \psi_j \, . &\nonumber
\end{alignat}

These equations imply some relations between the eigenfunctions of $\cA_G$ and those of $\cB$~: if $\eps_j$ is one of the two roots of
 $q\eps^2-\lambda_j \eps+1=0$, then $\tau_+ \psi_j-\eps_j \tau_- \psi_j$ is an eigenfunction of $\cB$ for the eigenvalue $\eps_j^{-1}$. It is possible to write a proof of Theorem \ref{t:thm1} using the eigenfunctions of $\cB$ instead of those of $\cA_G$, and this is the approach we will develop in \S \ref{sec:schro} for general graphs.
 \end{rem}

\subsection{Conclusion}
 
The first proof of quantum ergodicity for the adjacency matrix on regular graphs used microlocal analysis on trees \cite{ALM}. The proof we presented first appeared in \cite{A}. 
We see the advantage of working with the non-backtracking random walk $\cB$~: the exact orthogonality in \eqref{eq:ortho} comes from the fact that the operators $(\cM^{\ast\,r}K)_B= q^{-r} \cB^r K_B \cB^r$ live in orthogonal spaces $\mathscr{H}_j$ for different values of $r$.

The proof given here bears similarities with a proof by Brooks, Lindenstrauss and Le Masson \cite{BLL} in the special case where $R=0$, i.e. $K$ is a function on the vertex set $V_N$. The main idea there is the following : it is true that $\frac{1}{\lambda_j^r}\cA_G^r$ is difficult to analyze, but perhaps $P_r(\cA_G)$ is more approachable for a good polynomial $P_r$. The authors of \cite{BLL} choose the Chebyshev polynomial $P_r(\cos\theta) = \cos(r\theta)$. They start by observing that for any $\theta\in [0,\pi]$ and $n\ge 10$, one has $|\frac{1}{n} \sum_{r=1}^n\cos(r\theta)^2| \ge 0.3$. From this, they deduce that $\var(K) \lesssim \frac{100}{9}\|\frac{1}{n}\sum_{r=1}^n P_{2r}(\frac{\cA}{2\sqrt{q}})K P_{2r}(\frac{\cA}{2\sqrt{q}})\|_{\mathscr{H}}^2$. The various terms $P_{2r}(\frac{\cA}{2\sqrt{q}})K P_{2r}(\frac{\cA}{2\sqrt{q}})$, for $r\not= r'$, are not exactly orthogonal, but by properly rearranging the sum, some orthogonality property takes place and the theorem follows.


A variant would be to replace the Chebyshev polynomial $P_r$ by the spherical function $\Phi(\cdot, r)$ \eqref{e:spheri} -- so that $\Phi(\cA_G, r)$ is a renormalized  averaging operator over spheres of radius $r$ -- or by $\Psi_r$, where $\Psi_r$ is chosen so that $\Psi_r(\cA_G)$ is an averaging operator over balls of radius $r$. This is the path followed by Le Masson--Sahlsten \cite{LMS} to prove a quantum ergodicity result similar to Theorem \ref{t:thm1} for sequences of hyperbolic surfaces with genus going to $+\infty$.

\section{The general case}\label{sec:schro}
\subsection{The result}
We now turn to the general case of Schr\"odinger operators on arbitrary graphs of bounded degrees. Let $G_N=(V_N,E_N)$ be a graph with $|V_N|=N$ vertices and degree bounded by $D$. Consider a Schr\"odinger operator of the form $H_{G_N} = \cA_{G_N} + W_N$ where $\cA_{G_N} $ is the adjacency operator, and $W_N : V_N\to \R$ is a real-valued ``potential''. For simplicity, we assume that there exists $A$ such that $W_N(x)\in [-A,A]$ for all $x$, for all $N$. It will be convenient to regard $(G_N,W_N)$ as a ``coloured graph'', the map $W_N:V_N\to [-A,A]$ is viewed as a colour.

The theory of local weak convergence extends to this framework of coloured graphs without difficulty \cite{AL}. 
By choosing a root $x\in V_N$ uniformly at random, the deterministic coloured graph $(G_N, W_N)$ is turned into a random rooted coloured graph $(G_N, x, W_N)$.
The local weak convergence of $(G_N, W_N)$ means that the random rooted coloured graph $(G_N, x, W_N)$ converges in law. We denote by $[\cT,o,\cW]$ the limiting random variable, and $\prob$ its law. Thus, $\cT$ is a random graph with a random root $o\in\cT$ and a ``random potential'' $\cW : \cT\To [-A, A]$. In particular, the definition implies that the value distribution of $W_N$ converges to the law of the random variable $\cW(o)$.

Our general result may be roughly summarized as follows~: if the ``limit'' Schr\"odinger operator has purely absolutely-continuous spectrum ($\prob$-almost surely), then quantum ergodicity holds for any sequence converging to it.

This holds under some additional assumptions to be detailed below.

First, we assume that $\cT$ is $\prob$-almost surely a {\em{tree}}. This is equivalent to assumption \textbf{(BST)} in the previous section, namely that $G_N$ has ``few short loops''. In case of $(q+1)$-regular graphs, assumption \textbf{(BST)} implied the limit measure $\prob$ was concentrated on a single tree $\T_q$, here we assume more generally that $\prob$ is concentrated on the set of coloured rooted trees.

Next, we assume we have absolutely continuous (AC) spectrum at the limit. More specifically, we assume the imaginary part of the limit Green function has inverse moments. To state this precisely, we introduce some notation.

Given a coloured tree $[\cT,o,\cW]$, define the Schr\"odinger operator $\cH = \cA + \cW$ on $\ell^2(\cT)$. Note that it is random if $[\cT,o,\cW]$ is so. We denote its Green function by $\cG^{\gamma}(v,w) = \langle \delta_v,(\cH-\gamma)^{-1}\delta_w\rangle_{\ell^2(\cT)}$. For $v,w\in\cT$ such that $v\sim w$, denote by $\cT^{(v|w)}$ the tree obtained from $\cT$ by removing the branch emerging from $v$ that passes through $w$. Let $\cH^{(v|w)}$ be the restricted operator $\cH^{(v|w)}(u,u')=\cH(u,u')$ if $u,u'\in\cT^{(v|w)}$ and zero otherwise. The corresponding Green function is denoted by $\cG^{(v|w)}(\cdot,\cdot;\gamma)$. We then put $\hat{\zeta}_w^{\gamma}(v) = -\cG^{(v|w)}(v,v;\gamma)$ and assume~:

\medskip

There is a non-empty open set $I_1$ such that for all $s>0$,
\begin{equation}\label{e:green}
\sup_{\lambda\in I_1,\eta_0\in (0,1)} \expect\left(\sum_{o':\,o'\sim o} |\Im \hat{\zeta}_o^{\lambda+i\eta_0}(o')|^{-s}\right) <\infty \,.
\end{equation}

\medskip

We refer to this condition as \textbf{(Green)}. Here, $\expect$ is integration with respect to $\prob$, so we integrate over the coloured rooted trees $[\cT,o,\cW]$. It can be shown that condition \textbf{(Green)} implies that for $\prob$-a.e. $[\cT,o,\cW]$, the limiting Schr\"odinger operator $\cH$ has pure AC spectrum in $I_1$. 

\begin{rem} \label{r:reg}In the case of regular graphs treated in the previous section, $[\cT,o]$ is deterministic, equal to the regular tree $[\T_q, o]$ with an arbitrary root, and $\cW=0$. For $\gamma\not\in \R$,
it can be shown that $\hat{\zeta}_w^{\gamma}(v)$ is independent of $(v, w)$, equal to one of the roots of
$$q\zeta^2 -\gamma \zeta+1=0,$$
more precisely, the root with negative imaginary part if $\Im \gamma>0$ (see \eqref{e:zeta}). For $\gamma$ approaching the real axis, we see that the two roots become real if $|\gamma|\geq 2\sqrt{q}$, and purely imaginary if
$|\gamma|< 2\sqrt{q}$. Condition \eqref{e:green} will hold where $\Im \zeta^\gamma$ stays away from zero, that is, for any open interval $I_1$ such that $\bar{I_1}\subset (-2\sqrt{q}, 2\sqrt{q})$.
\end{rem}

As in Section \ref{sec:adj}, we assume we have a sequence of expanders. For irregular graphs, this is formulated as the existence of a spectral gap for the generator of the simple random walk~: define $P_N:\C^{V_N}\to\C^{V_N}$ by $(P_N\psi)(x) = \frac{1}{d_N(x)}\sum_{y\sim x} \psi(y)$, where $d_N(x)$ is the degree of $x\in V_N$. Then $P_N$ is self-adjoint in $\ell^2(V_N,d_N)$ (the reference measure now assigns weight $d_N(x)$ to a vertex $x$). Being expanders means that the eigenvalue $1$ of $P_N$ is simple (so $G_N$ is connected) and that the spectrum of $P_N$ in $\ell^2(V_N,d_N)$ is contained in $[-1+\beta,1-\beta]\cup \{1\}$. We refer to this assumption as \textbf{(EXP)}.

We may now state the result. Let $(G_N,W_N)$ be a sequence of finite coloured graphs as above, and assume that conditions \textbf{(BST)}, \textbf{(Green)} and \textbf{(EXP)} are satisfied.


\begin{thm} \label{t:thm2}
Let $(\psi_j^{(N)})_{j=1}^N$ be an orthonormal basis of eigenfunctions of $H_{G_N}$ for $\ell^2(V_N)$ with eigenvalues $(\lambda_j^{(N)})_{j=1}^N$.
 Fix $R\in\N$ and let ${\mathbf K}_N:V_N\times V_N\to\C$ satisfy ${\mathbf K}_N(x,y)=0$ if $d(x,y)>R$ and $\sup_N\sup_{x,y\in V_N}|{\mathbf K}_N(x,y)|\le 1$. 

Then\footnote{Technically, the theorem holds if $\psi_j^{(N)}$ are real-valued. This assumption is not necessary if $R=0$. If $R>0$, we need more precisely that $\overline{\psi_j^{(N)}(x)}\psi_j^{(N)}(y)\in \R$ for any $(x,y)\in B$.} for any interval $I$ with ${I}\subset I_1$,
\begin{equation}\label{thm:schro}
\lim_{\eta_0\downarrow 0} \lim_{N\to \infty} \frac{1}{N} \sum_{\lambda_j^{(N)}\in I} \left| \langle \psi_j^{(N)},{\mathbf K}_N\psi_j^{(N)}\rangle - \langle{\mathbf K}_N\rangle_{\lambda_j^{(N)}+i\eta_0}\right| =0\,,
\end{equation}
where $\langle {\mathbf K}\rangle_{\gamma} = \sum_{x,y\in V_N} {\mathbf K}(x,y)\Phi_{\gamma}^N({x},{y})$.
\end{thm}


As detailed below, the weight function $\Phi_{\gamma}^N(\tilde{x},\tilde{y})$ is expressed in terms of the Green function on the universal cover $\widetilde{G}_N$ of $G_N$. It is well-defined for $\gamma$ in the upper-half plane, which is the reason for the presence of 
the additional limit $\lim_{\eta_0\downarrow 0}$ (on a regular graph, one can directly take $\eta_0=0$, as the Green function on $\T_q$ stays finite on the real axis).  

Let $\widetilde{G}_N=(\widetilde{V}_N,\widetilde{E}_N)$ be the universal cover of $G_N$ (so $\widetilde{G}_N=\T_q$ in case the graphs are $(q+1)$-regular). Let $\widetilde{H}_N = \widetilde{\cA}_N + \widetilde{W}_N$ on $\ell^2(\widetilde{V}_N)$, where $\widetilde{\cA}_N$ is the adjacency matrix of $\widetilde{G}_N$, and $\widetilde{W}_N$ is the lift of $W_N$. We denote the corresponding Green function by $\tilg_N^{\gamma}(x,y) = \langle \delta_x,(\widetilde{H}_N-\gamma)^{-1}\delta_y\rangle_{\ell^2(\widetilde{V}_N)}$. Then
\[
\Phi_{\gamma}^N({x},{y}) = \frac{\Im \tilg_N^{\gamma}(\tilde{x},\tilde{y})}{\sum_{x\in V_N} \Im \tilg_N^{\gamma}(\tilde{x},\tilde{x})} \,,
\]
where $\tilde{x},\tilde{y}\in\widetilde{V}_N$ are lifts of $x,y\in V_N$ satisfying $d_{\widetilde{G}_N}(\tilde{x},\tilde{y}) =d_{G_N}(x,y) $.


In the very special case of an adjacency matrix over regular graphs, this reduces to $\Phi_{\gamma}^N({x},{y}) = \frac{1}{N} \frac{\Im \cG^{\gamma}(\tilde{x},\tilde{y})}{\Im \cG^{\gamma}(\tilde{x},\tilde{x})}$, and this quotient of Green functions is just another expression for the spherical function \eqref{e:spheri}. Note that for regular graphs, $\Phi_{\gamma}^N({x},{x})= \frac{1}{N}$ is the uniform measure; as we have seen, the result implies that for any given $a_N : V_N\To \C$, $\sum_{x\in V_N} a_N(x)| \psi_j^{(N)}(x)|^2$ is close to the uniform average $\frac{1}{N}\sum_{x\in V_N} a_N(x)$ for most $j$. For non-regular graphs, or in the presence of a potential $W_N$, it no longer holds true that $\Phi_{\gamma}^N({x},{x})= \frac{1}{N}$.
But for $a_N = \chi_{\Lambda_N}$, the characteristic function of a set $\Lambda_N\subset V_N$ of size $\geq \alpha N$, $0<\alpha\leq 1$,  we can nevertheless show that
$ \sum_{x\in \Lambda_N}\Phi_{\gamma}^N({x},{x})$ will always be bounded below by some $c_\alpha>0$,
with $c_\alpha$ depending only on $\alpha$. This implies that for most $j$, $\sum_{x\in  \Lambda_N}| \psi_j^{(N)}(x)|^2$ is also bounded below by $c_\alpha/2$.
 So our result can truly be interpreted as a delocalization result of eigenfunctions, saying that the mass of $| \psi_j^{(N)}(x)|^2$ on a macroscopic set is bounded below.



  \subsection{Discussion of assumptions}
 An assumption about the spectral gap (equivalently, the strong connectedness of the graphs), such as \textbf{(EXP)}, seems a reasonable replacement for the ergodicity assumption in the Shnirelman theorem. Obviously, quantum ergodicity does not hold for disconnected graphs.
On the other hand, assuming that the spectral gap is fixed is certainly too strong, a careful study of the proof reveals that one may allow the spectral gap to decay very slowly. Concrete examples reveal additional subtleties that would deserve further investigation~:
 discrete tori have a spectral gap decaying quite quickly
(polynomially with $N$). For these, eigenfunctions are explicit (trigonometric polynomials) and one can prove that quantum ergodicity holds if the test operator ${\mathbf K}_N$ is a multiplication operator; however, it does not hold in the form stated above if ${\mathbf K}_N$ is a derivation operator. 
  
  As we have said, \textbf{(Green)} implies, in a very strong manner, that the spectrum of the infinite Schr\"odinger operator $\cH = \cA + \cW$ on $\ell^2(\cT)$ is purely AC in the interval $I_1$. This is very reasonable, as AC spectrum is usually interpreted as meaning ``delocalization''. Our result is indeed a delocalization result~: it says that, if a sequence of finite systems converges to an infinite one having purely AC spectrum, then the eigenfunctions of the finite systems are delocalized, in the sense that $|\psi|^2$ is comparable to the uniform measure.
  
 On a more technical level, we actually need \textbf{(Green)} only for all $0< s \le s_0$, for some finite $s_0$ which in principle could be made explicit.
In addition, we believe that our results continue to hold under the following weaker variant of \textbf{(Green)}~: for all $s>0$, we have $\sup_{\eta_0\in (0,1)} \int_{I_1} \expect(\sum_{o'\sim o} |\Im \hat{\zeta}_o^{\lambda+i\eta_0}(o')|^{-s})\,\dd\lambda<\infty$.

Assumption \textbf{(BST)}, according to which our graphs have few short loops, seems to be of a merely technical nature, and we would like to work without it. However it is crucial in estimates such as
\eqref{e:2norms}, to show that the last term is negligible.

Finally, we mention that both assumptions \textbf{(EXP)} and \textbf{(BST)} are ``generic'' if $G_N$ is an $N$-lift of $G_0$. See \S~\ref{sec:reste} for details.

\subsection{The proof}

If we try to mimic the proof we gave in \S~\ref{sec:proof1} for the adjacency matrix on a regular graph, we note that for an irregular graph, $q$ becomes a function $q(x) = d(x)-1$, and the degrees mismatch, so that \eqref{e:smuggle} is no longer true. Analogous problems arise if the graph is regular but $W_N\neq 0$.
As noted in Remarks \ref{r:1} and \ref{r:2}, in the proof of the regular case, we indirectly considered powers of the non-backtracking operator \eqref{e:nonba}. We now build on this idea and try to convert the eigenfunctions $(\psi_j^{(N)})_{j=1}^N$ into eigenfunctions of some non-backtracking operator. The idea of considering a non-backtracking quantum variance first appeared in \cite{A}. Note that the idea of replacing simple random walks by non-backtracking ones is also useful to solve many other problems \cite{Fri, Bor-new, BLM, ABLS, LubPer}.

Denote $(\tau_-\varphi)(x_0,x_1)=\varphi(x_0)$ and $(\tau_+\varphi)(x_0,x_1)=\varphi(x_1)$. Let us try to find a transformation $\psi_j\rightsquigarrow f_j=\tau_+ \psi_j-\zeta_j \tau_- \psi_j$ (where $\zeta_j$ is now a function on $B$ instead of being a constant) such that 
$\cB f_j$ is proportional to $f_j$. While this appears to be impossible if we require $\cB f_j/f_j$ to be a constant function, we find that if $\zeta_j$ satisfies the system of $|B|$ algebraic equations
\begin{equation}\label{e:zetaj}
\lambda_j ={W}_N(v) + \sum_{u\in \cN_v\setminus\{w\}}\zeta(v, u) + \frac{1}{\zeta(w, v)} \,
\end{equation}
for all $(w, v)\in B$ then $\cB f_j = \zeta_j^{-1} f_j$. Thus we want to examinate the existence of solutions for these algebraic equations, but also their behaviour as $N\To +\infty$.

We recognize in \eqref{e:zetaj} the resolvent identity for the Schr\"odinger operator $\widetilde{H}_N=\widetilde{\cA}_N+\widetilde{W}_N$ on the universal cover $\widetilde G_N$ of $G_N$. More precisely, for $\gamma\in \C\setminus \R$, let
$\zeta^{\gamma}(y,x) = - \langle \delta_{\tilde{x}}, (\widetilde{H}_N^{(\tilde{x}|\tilde{y})}-\gamma)^{-1}\delta_{\tilde{x}}\rangle$.
As is well-known, it satisfies the recursion relation
\begin{equation}\label{e:zeta}
\gamma =  \widetilde{W}_N(v) + \sum_{u\in \cN_v\setminus\{w\}}\zeta^{\gamma}(v, u) + \frac{1}{\zeta^{\gamma}(w, v)} \,
\end{equation}
which is exactly the one we need in \eqref{e:zetaj}. As we do not know if $\zeta^{\gamma}(y,x)$ has a finite limit as $\gamma$ approaches the real axis, we let $\gamma_j=\lambda_j+i\eta_0$ for some fixed (arbitrarily small) $\eta_0>0$.
We henceforth denote $\zeta_y^{\gamma}(x) :=\zeta^{\gamma}(y,x)$ and let
\begin{equation}\label{e:fpsi}
f_j(x_0,x_1) =\psi_j(x_1) -  \zeta^{\gamma_j}_{x_0}(x_1)\psi_j(x_0) \,, \qquad f_j^{\ast}(x_0,x_1) = \psi_j(x_0) - \zeta^{\gamma_j}_{x_1}(x_0)\psi_j(x_1) \,.
\end{equation}
These functions satisfy $\zeta^{\gamma_j}\cB f_j =f_j - i\eta_0\tau_+\psi_j$ and $\iota \zeta^{\gamma_j}\cB^{\ast}f_j^{\ast} =  f_j^{\ast} -i\eta_0 \tau_-\psi_j$, where $\iota\zeta^{\gamma}(x_0,x_1) = \zeta^{\gamma}_{x_1} (x_0)$. As we let $\eta_0\downarrow 0$ at the end of the proof, the terms proportional to $\eta_0$ may be considered as negligible, and we omit them in this sketch of proof (in the regular case, we can take $\eta_0=0$ from the start, since $\zeta^\gamma$ is known to have a finite limit when $\gamma$ approaches the real axis, see Remark \ref{r:reg}).

Unlike the case of regular graphs, the functions $\zeta^{\gamma_j}$ now depend on $N$ (although this is absent from our notation) and we need to control their behaviour as $N\To +\infty$. More precisely we need to control the large and small values. Under assumption \textbf{(BST)}, the values of  $\zeta^{\gamma}$ converge {\emph{in distribution}} to
$\hat{\zeta}^\gamma$. The moment assumption \textbf{(Green)} \eqref{e:green} is precisely what will allow us to control the extreme values of $\zeta^{\gamma}$.

We are now ready to define the \emph{``non-backtracking'' quantum variance}, defined directly in terms of the functions $f_j$. Like in \S \ref{s:op}, instead of matrices we consider elements $K$ of $\mathscr{H}_k$, and we let\footnote{The word ``variance'' is perhaps unfortunate here as we have removed the square from the terms in the sum for technical reasons.}
\[
\varnbi(K) = \frac{1}{N} \sum_{\lambda_j\in I} \left|\langle f_j^{\ast},K_B f_j\rangle\right|,
\]
where $K_B$ is defined in \eqref{e:kb}. As in \S~\ref{sec:proof1}, we use the eigenfunction equation to show that the quantum variance is invariant under certain transformations. Namely, neglecting errors proportional to $\eta_0$, we have $\zeta^{\gamma_j}\cB f_j =f_j$ and $\iota \zeta^{\gamma_j}\cB^{\ast}f_j^{\ast} =  f_j^{\ast}$. This implies that
\[
\langle f_j^{\ast},K_Bf_j\rangle = \langle f_j^{\ast},(\cB \overline{ \iota \zeta^{\gamma_j}})^k K_B ( \zeta^{\gamma_j}\cB)^\ell f_j\rangle
\]
for any integers $k,\ell$, which replaces \eqref{e:smuggle}. Hoping as in \eqref{e:invariantrel} to have decay of correlations between the various terms, we write for any integer $T$,
\[
\varnbi(K) = \frac{1}{N} \sum_{\lambda_j\in I} \left|\left\langle f_j^{\ast},\frac{1}T\sum_{k=0}^{T-1}(\cB\overline{ \iota \zeta^{\gamma_j}})^{T-1-k} K_B ( \zeta^{\gamma_j}\cB)^k f_j\right\rangle\right|.
\]
The rest of the proof is substantially more involved than the one in \S~\ref{sec:proof1}. We roughly describe the main ideas.

The upper bound in terms of a Hilbert-Schmidt norm was previously very simple~: we combined the easy bound $\var(K)\le \|K_G\|_{HSN}^2$ with inequality \eqref{e:2norms}.   
Here, we have two difficulties. First, the family $(f_j)$ is not an orthonormal basis of $\ell^2(B)$.
In fact, $\|f_j\|$ is not even necessarily bounded for each $j$ as $\eta_0\downarrow 0$ (nevertheless, our assumptions imply that the mean $\frac{1}{N}\sum_{\lambda_j\in I}\|f_j\|^2$ stays bounded).
The second difficulty is the fact that $(\cB\overline{\iota \zeta^{\gamma_j}})^{T-1-k} K_B ( \zeta^{\gamma_j}\cB)^k$ depends on the eigenvalue $\lambda_j$. As a starting point to face these problems, we use the holomorphic functional calculus of operators to extend the operator $f_j\mapsto (\cB\overline{ \iota\zeta^{\gamma_j}})^{T-1-k} K_B ( \zeta^{\gamma_j}\cB)^k f_j$ to the whole of $\ell^2(B)$. 

Omitting technicalities, we are able to obtain an upper bound of the form
\[
\varnbi(K)^2 \lesssim \int_{\Re \gamma \in I, \Im \gamma=\eta_0} \left\| \frac{1}{T} \sum_{k=0}^{T-1} (\cB\overline{ \iota\zeta^{\gamma}})^{T-1-k} K_B ( \zeta^{\gamma}\cB)^k\right\|_{\gamma}^2\,\dd\gamma
\]
for some $\gamma$-dependent HS norm $\|\cdot\|_{\gamma}$. While in \S~\ref{sec:proof1}, the terms in \eqref{eq:ortho} were orthogonal to each other and we used Pythagoras, this won't be the case here, but the contribution of the off-diagonal terms (with $k'\neq k$) will decay exponentially in $|k-k'|$.

Decomposing the norm, we find that
\begin{multline}\label{e:integral}
\varnbi(K)^2 \lesssim \int_{\Re \gamma \in I, \Im \gamma=\eta_0} \frac{1}{T^2}\bigg(\sum_{k'\le k \le T-1}\langle (\cS_{\gamma}e^{i\theta_\gamma})^{k-k'} C_\gamma K, C_{\gamma}K\rangle_{\ell^2(B_k, \nu_\gamma)} \\ + \sum_{k<k'\le T-1} \langle C_{\gamma}K,(\cS_{\gamma}e^{i\theta_{\gamma}})^{k'-k}C_{\gamma}K\rangle_{\ell^2(B_k, \nu_\gamma)}\bigg) \dd\gamma\,,
\end{multline}
where $\cS_\gamma$ is the transfer operator $\C^{B_k}\To \C^{B_k}$, defined by
\[
(\cS_{\gamma}K)(x_0;x_k) = \frac{|\zeta^{\gamma}_{x_1}( x_0)|^2}{|\Im \zeta^{\gamma}_{x_1}( x_0)|} \sum_{x_{-1}\in \cN_{x_0}\setminus \{x_1\}} |\Im \zeta^{\gamma}_{x_0}( x_{-1})| \,K(x_{-1};x_{k-1})\, 
\]
and the function $e^{i\theta_\gamma}$ is a function of modulus $1$, defined by
\[
e^{i\theta_\gamma}(x_0,\ldots, x_k)=\overline{\zeta^\gamma_{x_0}( x_1)}^{-1}\zeta^\gamma_{x_0}( x_1)\,.
\]
Finally $C_\gamma(x_0,\ldots, x_k)=-\frac{1}{2\tilde g_N^\gamma(x_0, x_k)}$ is the inverse of the Green function on $\widetilde{G}_N$.

The operator $\cS_{\gamma}$ is sub-stochastic; in fact, it is exactly stochastic if we take $\eta_0=0$, as follows from the
relation
 \begin{equation}\label{e:sumzeta}
\sum_{u\in\cN_v\setminus\{w\}} |\Im \zeta^{\gamma}_v( u)| = \frac{|\Im \zeta^{\gamma}_w( v)|}{|\zeta^{\gamma}_w( v)|^2} - \eta_0\,
\end{equation}
which is a direct consequence of \eqref{e:zeta}. In \eqref{e:integral}, the measure $\nu_\gamma$ is the invariant probability measure for the Markov chain on $B_k$ defined by $\cS_{\gamma}$.

It is known by Wielandt's theorem that adding phases to a matrix with positive entries strictly diminishes its spectral radius, unless $\theta_\gamma$ is cohomologous to a constant; so it should be expected that $\cS_{\gamma}e^{i\theta_\gamma}$ decorrelates faster than $\cS_{\gamma}$. We do not apply directly Wielandt's theorem, because we are interested in operator norms, instead of spectral radii. Roughly speaking, using assumption \textbf{(Green)} to control the extreme values of $\zeta^{\gamma}$, and \textbf{(EXP)} to control spectral gaps, what we prove is that
\begin{itemize}
\item $\cS_{\gamma}$ has norm $1$ (as it is stochastic) and has a spectral gap on the orthogonal of constants 
\item typically, the presence of the phase $e^{i\theta_\gamma}$ will strictly diminish the norm of the fourth power $(\cS_{\gamma}e^{i\theta_\gamma})^4$, so each inner product in \eqref{e:integral} decays exponentially in $k-k'$ 
\item for non-typical situations, we do not control individual terms, but show instead that the phase $e^{i\theta_\gamma}$ induces cancellations between those terms, so that the mean sum of inner products decays with $T$, except when the potential $\cW=E_0$ is deterministic and $\gamma$ approaches $E_0$ (in that case, we discard that isolated value from our analysis).
\end{itemize}
Hence in any case, we finally establish a bound of the form
\begin{equation}\label{e:fini}
\varnbi(K)^2 \lesssim  \frac{C(\beta)}{T^{\alpha}} \left\Vert K\right\Vert^2_\infty +o_T(1)_{N\To +\infty} \left\Vert K\right\Vert^2_\infty
\end{equation}
for some $\alpha\in(0,1)$.


The proof is almost complete. Summarizing, we have shown that $\varnbi(K) \to 0$ as $N\to \infty$, for any $K\in\mathscr{H}_k$. What remains to be done is to connect this with the main result \eqref{thm:schro} via the formulas \eqref{e:fpsi} that relate the original eigenfunctions to the non-backtracking ones. For instance, in the special case of the adjacency matrix on regular graphs, our argument in this step is essentially to show that if $\varnbi(K) \to 0$, then $\vari(\nabla^{\ast}K)\to 0$, which implies the theorem by \eqref{e:red}.

One may wonder if there is a more direct proof of quantum ergodicity for non-regular graphs, which does not involve the non-backtracking quantum variance, but a generalized version of $\vari(\nabla^{\ast}K)$ instead.

\subsection{Conclusion}

We have shown that if a (random) Schr\"odinger operator $\cH$ on a (random) tree $\cT$ has pure AC spectrum, then his eigenvectors are delocalized, in the sense that if we consider a sequence of finite graphs $(G_N)$ converging to it, then the eigenfunctions of the finite model become uniformly distributed as $N\to\infty$ (in a weak sense).

A different approach would be to study directly the generalized eigenfunctions of $\cH$ on $\cT$. If the spectrum is purely AC, then the corresponding eigenfunctions are not in $\ell^2(\cT)$. This is a weak, immediate delocalization property for such eigenfunctions. Can one obtain more ? That is, do the generalized eigenfunctions of $\cH$ on $\cT$ inherit more specific delocalization properties, due to the fact that the eigenfunctions on the converging finite graphs are quantum ergodic ?

One should not imagine that all the generalized eigenfunctions of $\cH$ will be uniformly distributed on $\cT$. In fact consider $\cH=\cA_{\T_q}$ on $\T_q$, for which we proved a quantum ergodicity theorem in Section~\ref{sec:adj}. Fix any $o\in\T_q$ and consider the Poisson kernel $P_{\gamma,\xi}$ on $\T_q$ defined as follows~: given an infinite path $\xi=(o,v_1,v_2,\dots)$ in $\T_q$, let $P_{\gamma,\xi}(v) = \frac{\cG^{\gamma}(v\wedge \xi,v)}{\cG^{\gamma}(o,v\wedge \xi)}$, where $v\wedge \xi$ is the vertex of maximal length in $[o,v]\cap[o,\xi]$. Here $[o,v]$ is the path from $o$ to $v$ and $[o,\xi]=(o,v_1,v_2,\dots)$. Then $P_{\gamma,\xi}$ is a generalized eigenfunction of $\cA_{\T_q}$. In fact, such eigenfunctions form a ``basis'' for the set of generalized eigenfunctions of $\cA_{\T_q}$ (cf. \cite{AS} and references therein for more general operators). However, as the Green function $\cG^{\gamma}$ of $\cA_{\T_q}$ is explicit, by simple calculation one sees that for any $\lambda\in(-2\sqrt{q},2\sqrt{q})$, $|P_{\lambda+i0,\xi}(v)|^2$ grows exponentially as $v$ moves on the path $(o,v_1,v_2,\dots)$, and decays exponentially as $v$ moves away from it. 

Still, one could try to understand certain eigenfunctions of $\cA_{\T_q}$ using coloured graphs $(G_N,\psi_{\lambda_n^{(N)}})$, where $G_N$ satisfies \textbf{(BST)}, and the ``colour'' $W_N$ is now replaced by an eigenfunction $\psi_{\lambda_n^{(N)}}$ of $\cA_{G_N}$ corresponding to an eigenvalue $\lambda_n^{(N)}$, perhaps in a spectral window $[\lambda_{\ast}-\eps_N,\lambda_{\ast}+\eps_N]$. One can then try to understand the Benjamini-Schramm convergence of such coloured graphs, the limit should intuitively feature a random eigenfunction $\psi$ of $\cA_{\T_q}$. Perhaps such random eigenfunctions are uniformly distributed on $\T_q$. Ideally, if $\psi_{\lambda_n^{(N)}}$ is normalized as $\|\psi_{\lambda_n^{(N)}}\|_{\ell^2(V_N)}^2 = N$, one would like to have $|\psi(v)|^2 \approx 1$ for all $v\in \T_q$. The study of this problem, and a proof even of a weaker result $\|\chi_{\Lambda} \psi\|_{\ell^2(\cT)}^2\approx |\Lambda|$ for operators $\cH$ in general seems quite interesting.

\section{Trees of finite cone type and their stochastic perturbations}\label{sec:cone}

So far we derived a very general ``black box'' result asserting that if a sequence of Schr\"odinger operators on finite graphs converges to
a Schr\"odinger operator on a tree with pure AC spectrum, then quantum ergodicity holds. We need to find examples other than the one already covered by Theorem \ref{t:thm1}.

\subsection{Anderson model on the regular tree}\label{sec:andreg}

Let $\T_q$ be the $(q+1)$-regular tree. Consider a random Schr\"odinger operator $\cH_{\eps}=\cA+\eps \cW$ on $\T_q$, where $\cW(v)$ are i.i.d random variables.

Under very weak assumptions on the law of $\cW(v)$ (existence of a second moment is enough), it was shown by Klein \cite{Klein} that, given $0<E_0< 2\sqrt{q}$, the spectrum of $\cH_{\eps}$ in $[-E_0, E_0]$ is a.s. purely AC, for all $\eps$ small enough. This can be strengthened to obtain condition \textbf{(Green)} if we put stronger assumptions on the law of $\cW(v)$~:
\medskip

\textbf{(POT)} The $(\cW(v))_{v\in \T_q}$ are i.i.d. with common distribution $\nu$ which has a compact support $\supp \nu \subseteq [-A,A]$, and is H\"older continuous, i.e. there exist $C_{\nu}>0$ and $b\in (0,1]$ such that for any bounded interval $J\subset \R$, $\nu(J) \le C_{\nu} \cdot |J|^b$.

\medskip

Under this assumption, an argument from \cite{AW2} shows that the \textbf{(Green)} bound \eqref{e:green} holds on $[-E_0, E_0]$, if $\expect$ is the expectation with respect to the random potential on $\T_q$ and $\eps>0$ is small enough. Thus, Theorem \ref{t:thm2} applies on the interval $I=[-E_0, E_0]$, for any sequence $(G_N,\eps W_N)$ of finite coloured graphs such that
\begin{itemize}
\item[--] the graphs $(G_N)$ satisfy assumptions \textbf{(BST)} and \textbf{(EXP)},
\item[--] $(G_N,\eps W_N)$ converges to $[\T_q,o,\eps\cW]$, in the sense of Benjamini-Schramm, 
\item[--] the law of $(\cW(v))_{v\in\T_q}$ satisfies \textbf{(POT)}, and $\eps$ is small.

\end{itemize}
For example if we choose $W_N=(W_N(x))_{x\in G_N}$ i.i.d. random variables of law $\nu$, which are independent for different $N$, it is shown in \cite{AS3} that the second condition is satisfied for almost every realization of $(W_N)$. Moreover, \emph{if $K_N$ is deterministic, or probabilistically independent of $W_N$}, the average $\langle K_N\rangle_{\gamma}$ can be replaced by a simplified expression. In particular, if $K_N=a_N$ is a function, $\langle a_N\rangle_{\gamma}$ can be replaced by $\langle a_N\rangle$, the uniform average.

\subsection{Trees of finite cone type}
In the previous subsection we used Theorem~\ref{t:thm2} to prove an aspect of spatial delocalization for the Anderson model on $\T_q$. We now show that Theorem~\ref{t:thm2} can also be used to study delocalization on some irregular trees $\cT$. We focus here on the simpler case of the adjacency matrix on $\cT$. The Anderson model on irregular trees is deferred to \S~\ref{sec:Andersoncone}.

Trees \emph{of finite cone type} \cite{Nag, KLW3} are rooted trees $(\cT, o)$ satisfying the following condition. Given $v\in \cT$, define the cone $\mathscr{C}(v) = \{w\in \cT : v\in [o,w]\}$, where $[o,w]$ is the unique path from $o$ to $w$. So $\mathscr{C}(v)$ is the forward subtree emanating from $v$ (and $\mathscr{C}(o)=\cT$ is seen as a cone from $o$). We say that $(\cT, o)$ is of finite cone type if the number of non-isomorphic cones is finite. Such trees are sometimes called \emph{periodic trees} \cite{LP16}, and emerge when considering spanning trees of regular tessellations of the hyperbolic plane.

One can show that every finite directed graph $(G,x_0)$ has a cover $(\cT,x_0)$ which is a tree of finite cone type, and conversely, every tree of finite cone type $(\cT,o)$ covers a finite directed graph $(G,x_0)$. Universal covers of finite undirected graphs (sometimes called \emph{uniform trees} or {quasi-homogeneous trees}) are also of finite cone type; see \cite{Nag}.

We may state the following particular case of Theorem \ref{t:thm2}~:

\begin{cor}\label{cor:cone}
Let $(G_N)$ be a sequence of graphs of degree bounded by $D$, converging to a tree $\cT$ of finite cone type in the local weak sense. Assume the adjacency matrix on $\cT$ satisfies \emph{\textbf{(Green)}} on some open set $I$, and assume the $(G_N)$ satisfy \emph{\textbf{(EXP)}}. Then for any orthonormal basis $(\psi_j^{(N)})$ of eigenfunctions of $\cA_{G_N}$, quantum ergodicity holds in $I$, that is, \eqref{thm:schro} holds true.
\end{cor}

More precisely, we mean that $(G_N)$ has a Benjamini-Schramm limit $\prob$ which is concentrated on $\{[\cT,v]\}_{v\in \cT}$ (where $\cT$ is fixed and each $(\cT,v)$ is of finite cone type), and \textbf{(Green)} is considered with respect to this measure $\prob$.

To show that the previous result is non-void, we dedicate \S~\ref{s:Puiseux} and \S~\ref{sec:reste} to show that there exist trees $\cT$ of finite cone type (other than $\T_q$) satisfying \textbf{(Green)} in large regions of the spectrum, and sequences $(G_N)$ which satisfy \textbf{(EXP)} and converge to $\cT$ in the local weak sense.

\subsection{Assumption \textbf{(Green)} for trees of finite cone type \label{s:Puiseux}}

There are already some results in the literature which consider the problem of spectral delocalization on irregular trees. In \cite{Ao91}, conditions are given to exclude point spectrum in uniform trees. The paper \cite{KLW2} studies a different class of trees of finite cone type, and achieves a good control on the Green function of $\cA_{\cT}$, which suggests that assumption \textbf{(Green)} holds true through most of $\sigma(\cA_{\cT})$. However, it seems that the trees of \cite{KLW2} are never \emph{unimodular}\footnote{We say that a (fixed, deterministic) tree $\cT$ is unimodular if there exists a unimodular measure $\prob$ which is concentrated on $\{[\cT,v]\}_{v\in \cT}$. Here, $\prob$ is unimodular if it satisfies the ``mass-transport principle'' $\int \sum_{o'\sim o} f([G,o,o'])\,\dd\prob([G,o]) = \int\sum_{o'\sim o} f([G,o',o])\,\dd\prob([G,o])$. See \cite{AL} for details.}. If the reader is not familiar with this notion, we just mention that the assumptions of Corollary~\ref{cor:cone} imply that the tree $\cT$ must be unimodular. So our aim in what follows is to prove that assumption \textbf{(Green)} holds true on large regions of the spectrum of trees of finite cone type, under less restrictive assumptions than in the paper \cite{KLW2}. Concrete examples which fit the framework of Corollary~\ref{cor:cone} will then be provided.

Consider a finite set  of labels $\mathfrak{A} = \{1,\dots, m\}$ and a matrix $M=(M_{j,k})_{j,k\in \mathfrak{A}}$, where $M_{j,k}\in\N$. A tree $\T(M,j)$ is constructed by asserting that the root has the label $j$, and that each vertex with label $k$ has $M_{k,l}$ children of label $l$. Such trees have $m$ cone types, and any tree of finite cone type arises in this fashion.

We shall make the following assumption~:

\smallskip

\begin{enumerate}[\textbf{(C1)}]
\item We have $M_{1,1}=0$ (and $M_{1, k}>0$ for at least one $k$). Moreover, for any $k,l\in \{2,\dots, m\}$, there is $n=n(k,l)$ such that $(M^n)_{k,l}\ge 1$.
\end{enumerate}

\smallskip

We could also consider the variant where $M_{1,1}$ is arbitrary while the condition holds on the full set $k,l\in \{1,\dots, m\}$. This is called (M2) in \cite{KLW2}.

Assumption \textbf{(C1)} says that all cone types arise at some point as offspring of a given cone, except for the cone with label $1$ which plays a special role and may not reappear. In practice, this condition happens for the cone at the root, and has been introduced to allow for the example of the regular tree $\T_q$. It has two cone types : the cone at the origin has $(q+1)$ children while any other cone has $q$ children. Another such example is the $(p+1,q+1)$-biregular tree $\T_{p,q}$. More examples will be given later.

A tree of finite cone type $(\T,1)$ with associated $m\times m$ matrix $M=(M_{j,k})$ gives rise to the following system of polynomial equations
\begin{equation}\label{eq:pols}
\sum_{k=1}^m M_{j,k} h_k h_j -\gamma h_j + 1 = 0\,, \qquad j=1,\dots, m
\end{equation}
where $(h_1,\dots,h_m) \in \C^m$ and $\gamma\in\C$ is fixed; cf. \cite{KLW2}. If $\Im \gamma >0$, one of the solutions of the system is $(\zeta_1^\gamma,\dots,\zeta_m^\gamma)$, where $\zeta_1^{\gamma}=-G^{\gamma}(u,u)$ if $u$ has label $1$, and for $j\ge 2$, $\zeta_j^\gamma = \zeta_v^\gamma(v_+)$ if $v_+$ has label $j$. We now show that

\begin{prp}
\phantomsection\label{prp:Lang}
\begin{enumerate}[\rm (1)]
\item There is a discrete set $ \mathfrak{D}\subset \R$ such that, for all $j=1,\dots,m$, the solutions $h_j(\lambda+i\eta) = \zeta_j^{\lambda+i\eta}$ of \eqref{eq:pols} have a finite limit as $\eta \downarrow 0$ for all $\lambda\in \R\setminus  \mathfrak{D}$. The map $\lambda\mapsto \zeta_j^{\lambda+i0}$ is continuous on $\R\setminus \mathfrak{D}$, and there is a discrete set $\mathfrak{D}'$ such that it is analytic on $\R\setminus (\mathfrak{D}\cup\mathfrak{D}')$.
\item If moreover $(\T,1)$ satisfies \emph{\textbf{(C1)}}, then 
\begin{enumerate}[\rm(i)]
\item For any $j\ge 1$, the map $\lambda\in  \sigma(\cA_{\T})\setminus  \mathfrak{D}\mapsto | \Im\zeta_j^{\lambda+i0}|$ has finitely many zeroes.
\item $\sigma(\cA_{\T})$ is a finite union of closed intervals and points, $\cup_{r=1}^\ell I_r \cup \mathfrak{F}$. Moreover, the limit $\zeta_j^{\lambda+i0}$ exists on the interior $\mathring{I}_r$ and satisfies $|\Im \zeta_j^{\lambda+i0}|>0$ for all $j\ge 1$.
\end{enumerate}
\end{enumerate}
\end{prp}
  
As a consequence of part (2), for any $0<c<C<+\infty$, the set of $\lambda\in \sigma(\cA_{\T})$ such that $c< | \Im\zeta_j^{\lambda+i0}|<C$ for all $j$
is of the form $\sigma(\cA_{\T})\setminus I_{c, C}$ where $I_{c,C}$
is a finite union of open intervals, that shrinks to a finite set when $c\To 0$ and $C\To +\infty$. As we will see in the examples, (after defining the probability measure $\prob$) \textbf{(Green)} holds on $\sigma(\cA_{\T}) \setminus  I_{c, C}$.

\begin{proof}
To prove (1), it is sufficient to show that in a neighbourhood of any $\lambda_0\in \R$, the function $h_j(\gamma)=\zeta_j^{\gamma}$ has a convergent \emph{Puiseux expansion}, that is, an expansion of the form
\begin{equation}\label{eq:puiseux}
h_j(\gamma)=\sum_{n\geq m} a_n (\gamma-\lambda_0)^{n/d}
\end{equation}
for some $m\in \Z$ and some integer $d$, where the entire series $\sum_{n\geq 0} a_n z^n$ has a positive radius of convergence. The fact that $\zeta_j^{\gamma}$ possess such expansions is announced without explanation in a paper by Aomoto \cite{Ao91}; we give the full argument below.


Proposition VIII.5.3 in \cite{Lang} teaches us the following~: if $K$ is a field and if $(h_1,\ldots, h_m)$ is a solution, in some extension $K'$ of $K$, of a system of polynomial equations $P_1(h_1, \ldots, h_m)=0, \ldots, P_m(h_1, \ldots, h_m)=0$ where $P_j\in K[X_1, \ldots, X_m]$, and if the determinant of the Jacobian matrix $\det (\frac{\partial P_j}{\partial h_j}(h))\not=0$, then each $h_j$ is algebraic over $K$.

We apply this with
 $P_j(h_1,\dots,h_m) = \sum_{k=1}^m M_{j,k} h_k h_j -\gamma h_j + 1$. Clearly, we have $P_j \in K[h_1,\dots,h_m]$, where $K=\mathscr{K}_{\lambda_0}$ is the field of functions $f(\gamma)$ possessing a convergent Laurent series $f(\gamma)=\sum_{j=-n_0}^{\infty}a_j(\lambda_0) (\gamma-\lambda_0)^j$ in some neighbourhood $N_{\lambda_0}\subset \C$ of $\lambda_0$.

Let $K'=\mathscr{J}_{\lambda_0}$ be the field of functions $f$ which are meromorphic on $N_{\lambda_0}\cap \C^+$ for some neighbourhood $N_{\lambda_0}$ of $\lambda_0$. Then $K'$ is an extension of $K$, and we know that $\zeta_j^{\gamma}$ belong to $K'$ and satisfy $P_j(\zeta_1^{\gamma},\dots,\zeta_m^{\gamma})=0$. Calculating $J(h)=\det (\frac{\partial P_j}{\partial h_j}(h))$ and using that $\zeta_j^{\gamma}\sim \frac{1}{\gamma}$ as $\Im \gamma \to \infty$, we see that $J(\zeta^{\gamma}):=J(\zeta_1^{\gamma},\dots,\zeta_m^{\gamma})\neq 0$ for any $\gamma$ with a large enough imaginary part. Since $\gamma\mapsto J(\zeta^{\gamma})$ is holomorphic on $\C^+$, it follows that it cannot vanish identically on any neighbourhood $N_{\lambda_0}\cap \C^+$. Hence, $J(\zeta^{\gamma})$ is not the zero element of $K'$.

It follows that each $\zeta_j^{\gamma}$ is algebraic over $K$. By the Newton-Puiseux theorem (see e.g. \cite[Theorem 3.5.2]{Sim15}), each $\zeta_j^{\gamma}$ thus has an expansion of the form \eqref{eq:puiseux} in some neighbourhood $N_{\lambda_0}$ of $\lambda_0$. In particular, it is analytic near any $\lambda\in N_{\lambda_0}\setminus \{\lambda_0\}$. 
The set $ \mathfrak{D}$ corresponds to those $\lambda_0$ for which $m<0$ in the Newton-Puiseux expansion at $\lambda_0$, and the set $\mathfrak{D}'$ corresponds to those $\lambda_0$ for which $d>1$.


We now turn to (2). Denote $\zeta_j^{\lambda}:=\zeta_j^{\lambda+i0}$ for $\lambda\notin\mathfrak{D}$.

Given $\lambda_0\in \sigma(\cA)$, we expand $\lambda\mapsto\Im \zeta_j^{\lambda}$ in a disc $N_{\lambda_0}^{\eps}=\{|\gamma-\lambda_0|<\eps\}$ using the Newton-Puiseux expansion \eqref{eq:puiseux}. We take $\eps$ small enough so the expansion holds for all $j\ge 1$ and $\zeta_j^{\lambda}$ is well-defined except perhaps at $\lambda_0$. Note that $\zeta_j^{\lambda}\neq 0$, due to the recursive relation $\zeta_j^{\lambda} = \frac{-1}{\sum_{k=1}^mM_{j,k} \zeta_k^{\lambda}-\lambda}$.

Now $\Im \zeta_j^{\lambda} = \sum_{n\ge m} c_n(\lambda-\lambda_0)^{n/d}$, where $c_n=\Im a_n$ for $\lambda \in (\lambda_0,\lambda_0+\eps)$, and $c_n = \Im (a_n e^{i\pi n/d})$ for $\lambda \in (\lambda_0-\eps, \lambda_0)$. It follows that if $\Im \zeta_j^{\lambda}=0$ for some $\lambda\in (\lambda_0-\eps,\lambda_0+\eps)\setminus\{\lambda_0\}$, then one of the following situations must occur~:
\begin{itemize}
\item[--] either $\Im \zeta_j^{\lambda}=0$ on $O_{\lambda_0}=(\lambda_0-\eps, \lambda_0-\eps)\setminus\{\lambda_0\}$;
\item[--] or $\Im \zeta_j^{\lambda}=0$ on $O_{\lambda_0}'=(\lambda_0, \lambda_0+\eps)$ and is nonzero in some $(\lambda_0-\eps',\lambda_0)$;
\item[--] or $\Im \zeta_j^{\lambda}=0$ on $O_{\lambda_0}'=(\lambda_0-\eps, \lambda_0)$ and is nonzero in some $(\lambda_0,\lambda_0+\eps')$;
\item[--] or $\lambda_0$ is isolated from the zeroes of $\Im \zeta_j^{\lambda}$.
\end{itemize}

Suppose case 1 occurs for some $\Im \zeta_j^{\lambda}$, $j\ge 2$. Using \textbf{(C1)} and the argument in \cite[Lemma 4]{KLW2}, it follows that $\Im \zeta_k^{\lambda}=0$ on $O_{\lambda_0}$ for all $k\ge 1$. This in turn implies that $\Im G^{\lambda}(v,v)=0$ on $O_{\lambda_0}$ for all $v\in \T$ by \cite[Proposition 3]{KLW2}. Hence, either $\lambda_0\notin \sigma(\cA)$, or $\lambda_0$ is an isolated eigenvalue. If case 1 occurs for $\Im \zeta_1^{\lambda}$, we use $|\Im \zeta_1^{\lambda}| = (\sum_{k=2}^m M_{1,k}|\Im \zeta_k^{\lambda}|)|\zeta_1^{\lambda}|^2$ to reach the same conclusion.

The same analysis shows that if case 2 or 3 occurs for some $\Im \zeta_j^{\lambda}$, $j\ge 1$, then $\lambda_0$ must be a boundary point of $\sigma(\cA)$, which is isolated from the zeroes inside $\sigma(\cA)$.

We thus showed that the zeroes of $\Im \zeta_j^{\lambda}$, $j\ge 1$, do not accumulate near any $\lambda_0\in\sigma(\cA)$. This proves (2.i). For (2.ii), we use again the argument in Lemma 4 of \cite{KLW2} and their Proposition 3 to see that under \textbf{(C1)}, we have $\{\lambda\in \R\setminus \mathfrak{D} : \Im G^{\lambda}(v,v)>0 \} = \{\lambda\in\R\setminus\mathfrak{D} : |\Im \zeta_j^{\lambda}|>0\}$, where all sets coincide for arbitrary $v\in\cT$ and $j\ge 1$. In particular, $\sigma(\cA)\setminus \mathfrak{D} = \overline{\{\lambda\in \R\setminus \mathfrak{D} : \Im G^{\lambda}(v,v)>0\}}\setminus\mathfrak{D}$. On the other hand, $\{\lambda\in\R\setminus\mathfrak{D} : |\Im \zeta_j^{\lambda}|>0\}$ is a  finite union of intervals $\cup_{r=1}^\ell J_r$ by the continuity of $\R\setminus\mathfrak{D} \ni \lambda \mapsto \Im \zeta_j^\lambda$ and (2.i). Hence, if $\mathfrak{F} = \sigma(\cA)\cap \mathfrak{D}$, we get $\sigma(\cA) = (\sigma(\cA)\setminus \mathfrak{D})\cup \mathfrak{F}$, yielding (2.ii).
\end{proof}

\begin{rem}
In the previous proposition, we implicitly assumed that $\sigma(\cA_{\T})\setminus \mathfrak{D} \neq \emptyset$, i.e. the spectrum of the infinite tree $\T$ is not reduced to a finite set of points. This property holds in particular for the trees most relevant to us, namely universal covers of finite undirected graphs of minimal degree at least $2$. In this case, it follows from the results of \cite[\S 1.6]{BSV17} that $\sigma(\cA_{\T})$ has a continuous part.
\end{rem}

As mentioned in the beginning of this section, the assumptions of Corollary~\ref{cor:cone} imply that the tree $\cT$ must be unimodular. So let us conclude this section by giving examples of trees of this type.

Uniform trees (universal covers of finite undirected graphs) are unimodular  \cite[Example 9.3]{AL}. In fact, these appear to be the only unimodular trees of finite cone type \cite{Alt9}. As an example of a unimodular measure, one can take the uniform measure $\prob = \frac{1}{|G|}\sum_{x\in G} \delta_{[\widetilde{G},x]}$, where $G$ is the finite graph and $(\widetilde{G},x)$ its universal cover at $x\in G$. The unimodularity of this measure simply follows from the relation $\sum_{(x,y)\in B} f(x,y)= \sum_{(x,y)\in B} f(y,x)$, where $B$ are the directed edges of $G$.

To test for condition \textbf{(C1)}, it is not difficult to see \cite{Nag} that all cone types of such universal covers can be indexed by the directed edges of the finite graph $G$, except for the cone at the origin\footnote{This implies in particular that the number of cone types is bounded by $|B|+1$. The number of cone types is usually smaller than this, because many directed edges play an equivalent role in $G$. For example, if $G$ is regular, then $\widetilde{G}$ only has $1+1 = 2$ cone types, regardless of the size of $|B|$.}. Assumption \textbf{(C1)} says that any cone type, except perhaps for the cone at the origin, should arise as an offspring of a given cone. For this, it suffices to know that for any directed edges $b,b'$ in $G$, there exists a non-backtracking path $(v_0,\dots,v_k)$ in $G$ such that $(v_0,v_1)=b$ and $(v_{k-1},v_k)=b'$. In other words, it suffices to know that the non-backtracking operator is irreducible. This property is known to be true for any finite graph $G$ of minimal degree $\ge 2$ which is not a cycle, see \cite[Lemma 3.1]{OW07}.

To summarize, if $G$ is a finite connected graph of minimal degree $\ge 2$ and $G$ is not a cycle, then Proposition~\ref{prp:Lang} applies to $\widetilde{G}$, and \textbf{(Green)} holds true for $\prob = \frac{1}{|G|}\sum_{x\in G} \delta_{[\widetilde{G},x]}$.

Another way of generating unimodular trees without explicit reference to covers is as follows. Consider a matrix $A_{\cN} = \begin{bmatrix} a& b\\ c& d \end{bmatrix}$ with entries in $\N$. We construct a tree with two colours of vertices $\bullet,\circ$, according to the rule~: $\bullet$ has $a$ neighbours of type $\bullet$ and $b$ neighbours of type $\circ$, while $\circ$ has $c$ neighbours of type $\bullet$ and $d$ neighbours of type $\circ$. Such trees are unimodular if the root is chosen according to the probability measure $\prob = \frac{c}{b+c}\delta_{[\T,\bullet]} + \frac{b}{b+c} \delta_{[\T,\circ]}$. They can also be represented by cone matrices. For example, if all entries of $A_{\cN}$ are nonzero, then $(\T,\bullet)$ corresponds to the cone matrix $M =\begin{pmatrix} 0 & a & 0 & b & 0 \\ 0 & a' & 0 & b & 0 \\ 0 & a & 0 & b' & 0 \\ 0 & 0 & c' & 0 & d \\ 0 & 0 & c & 0 & d' \end{pmatrix}$, where $e'=e-1$ for $e=a,b,c,d$. Here, $\zeta_1^{\gamma} = -G^{\gamma}(\bullet,\bullet)$ while the $\zeta_j^{\gamma}$ for $j=2,\dots,5$ correspond to $\zeta_{\bullet}^{\gamma}(\bullet)$, $\zeta_{\circ}^{\gamma}(\bullet)$, $\zeta_{\bullet}^{\gamma}(\circ)$ and $\zeta_{\circ}^{\gamma}(\circ)$, respectively. Hence, \eqref{e:green} reads  $\expect(\sum_{y\sim o}|\Im \zeta_o^{\gamma}(y)|^{-s}) = \frac{c}{b+c}[a\,|\Im \zeta_{\bullet}^{\gamma}(\bullet)|^{-s}+b\,|\Im \zeta_{\bullet}^{\gamma}(\circ)|^{-s}] + \frac{b}{b+c}[c\,|\Im \zeta_{\circ}^{\gamma}(\bullet)|^{-s}+d\,|\Im \zeta_{\circ}^{\gamma}(\circ)|^{-s}]$. One checks that $M$ satisfies (C1) (with $n(k,l)=2$) if $\min(a,b,c,d)\ge 2$. So applying Proposition~\ref{prp:Lang}, we see \textbf{(Green)} holds on $\T$ endowed with $\prob$, on the set $\sigma(\cA_{\T})\setminus I_{c, C}$. One could also consider the same example with $d=0$ and $\min(a,b,c)\ge 2$ and obtain analogous results.

One may also use a larger number of colours. For example, with three colours, a $3\times 3$ matrix $A_{\cN}=[a_{ij}]$ with analogous rule gives rise to a unimodular tree as long as $a_{31}a_{23}a_{12}=a_{32}a_{13}a_{21}$. Note that if we construct an $n$-coloured tree following a ``neighbour matrix'' $A_{\cN}$ of size $n\times n$, then for any choice of the root the associated ``cone matrix'' $M$ appearing in \textbf{(C1)} is of size at most $(n^2+1)\times(n^2+1)$.

\subsection{About assumptions \textbf{(EXP)} and \textbf{(BST)}} \label{sec:reste}

It is known that any unimodular measure concentrated on (coloured) rooted trees is sofic \cite{E,BLS}. So once we have a random unimodular tree $\cT$, we may always assume there is a sequence $(G_N)$ which converges to $\cT$. In which cases can we ensure, in addition, that $(G_N)$ are expanders~?

Let $G_0$ be a given connected graph with minimal degree $d\ge 3$.
If $G_N$ is a random $N$-covering of $G_0$, it is known that asymptotically almost surely (a.a.s.), $(G_N)$ satisfies \textbf{(EXP)}. This is a standard combinatorial argument,
see \cite[\S 6.2]{ACF} for a counting argument of a same flavor in the case of random regular graphs. Substantially stronger arguments \cite{Fri03,Pu} give explicit bounds on the spectral gap of $G_N$ in terms of the spectral radius $\rho_P$ of the universal cover $\widetilde{G}_0$. Here $\rho_P$ is more precisely the spectral radius of the Laplacian $(Pf)(x) = \frac{1}{d(x)}\sum_{y\sim x}f(y)$ on $\ell^2(\widetilde{G}_0,d)$. A simple argument \cite{LN} yields $\rho_P \le \frac{2\sqrt{D-1}}{d}$, where $D$ and $d$ are the maximal and minimal degrees, respectively, which is optimal in the regular case where $d=D$. This combined with \cite{Fri03,Pu} gives an explicit bound on the spectral gap, holding a.a.s.

It is also easy to see that a random $N$-lift  has few short cycles a.a.s. and thus converges to the tree $\widetilde{G}_0$. More precisely, $(G_N)$ converges to the uniform measure $\prob = \frac{1}{|G_0|}\sum_{x\in G_0}\delta_{[\widetilde{G}_0,x]}$. See \cite[Lemma 24]{Bor-new} and \cite[Lemma 9]{BDGHT} for a proof.

Summarizing, we have shown that if $\T$ is any uniform tree satisfying \textbf{(C1)}, then by \S \ref{s:Puiseux} it satisfies  \textbf{(Green)} on the large set $\sigma(\cA_{\T})\setminus I_{c, C}$, and by \S \ref{sec:reste}, there exists a sequence $(G_N)$ satisfying both \textbf{(EXP)} and \textbf{(BST)} which converges to $\T$. Corollary~\ref{cor:cone} thus applies to this nontrivial context. As a bonus, we have shown that $G_N$ can be chosen as $N$-lifts of the finite graph which $\T$ covers, and that a random $N$-lift will do the job with high probability.

\subsection{Condition \textbf{(Green)} on random trees}\label{sec:Andersoncone}

As we mentioned in \S~\ref{sec:andreg} in case of the Anderson model on $\T_q$, if the random variables $(\cW(v))_{v\in\T_q}$ satisfy condition \textbf{(POT)}, then assumption \textbf{(Green)} holds true for $\cH_{\eps}=\cA+\eps \cW$ on any closed subset of $(-2\sqrt{q},2\sqrt{q})$, provided $\eps$ is small. Note that for $\eps=0$ \textbf{(Green)} holds precisely on $(-2\sqrt{q},2\sqrt{q})$. This suggests the following open question~

\begin{pro} \label{pro:pertu}
Is Assumption \textbf{(Green)} stable under small perturbations~? That is, if a coloured random rooted tree $(\cT, o, \cW)$ satisfies \textbf{(Green)} and we add a small random perturbation to $\cW$, will \textbf{(Green)} still hold~?
\end{pro}

Note that the spectrum of the Anderson model $\cH_{\eps}$ on $\Z$ is pure point for any $\eps\neq 0$, so the answer is negative in this case. It seems reasonable to assume the average degree of the root is larger than $2$ in order to get a positive answer.

One can also consider assumption \textbf{(Green)} for different models~:
\begin{pro} 
Does assumption \textbf{(Green)} hold for small enough $\eps$ in the following models ?
\begin{enumerate}[\rm (a)]
\item The Anderson model on $\T_q$ where the random potential only takes the values $0$ and $1$ with probability $1-\eps$ and $\eps$, respectively.
\item A Galton-Watson tree where the number of children is either $q_1$ or $q_2$ with probability $1-\eps$ and $\eps$, respectively, and the unimodular variant of this model.
\end{enumerate}
\end{pro}
In both models, we make a strong perturbation (varying $W$ from $0$ to $1$ or changing the tree configuration), but with small probability. In contrast, in \S~\ref{sec:andreg} we studied the effect of adding a uniformly small perturbation.

The argument given by Klein \cite{Klein} can be adapted to the first model to derive positive moment bounds on the Green function. It seems this argument can also be adapted to the second model to derive similar bounds; see also \cite{Kel12} for related results. However, we do not know if \textbf{(Green)} holds for these models, that is, if \eqref{e:green} also holds. Indeed, for the first model, the $(\cW(v))$ do not satisfy \textbf{(POT)}, since they follow the Bernoulli distribution $\nu_{\eps} = (1-\eps)\delta_0 + \eps \delta_1$.

\medskip

For the remaining of this section, we give a positive answer to Problem \ref{pro:pertu} for a special class of (unimodular) trees of finite cone type.

Let $\T$ be a tree of finite cone type with labels $\mathfrak{A}=\{1,\dots,m\}$ and cone matrix $M=(M_{j,k})$. Consider a Schr\"odinger operator $\cH_{\eps} =\cA+\eps \cW $ on $\T$, where $\cW (v)  $ are random variables satisfying Assumption \textbf{(POT)}.

We denote $\zeta_{j,\eps}^{\gamma} = \zeta_v^{\gamma}(v_+)$ if $v_+$ has label $j$. Though the function $\zeta_v^{\gamma}(v_+)$ depends on both vertices $v,v_+$, its distribution only depends on the type of $v_+$, which justifies the notation.

We assume that the minimal degree is at least $3$. We also denote $\mathbf{P}=\mathop\otimes_{v\in\T}\nu$ and $\mathbf{E}$ the corresponding expectation.

Fix $\delta\in (0,1)$. We introduce
\begin{equation}\label{e:sigmaac}
\sigma_{ac}^{\eps}(\delta) = \{\lambda \in \R : \mathbf{P}(|\Im \zeta_{j,\eps}^{\lambda+i\eta}|>\delta)>\delta \quad \forall j=1,\dots, m, \ \forall \eta \in (0,1)\} \, .
\end{equation}

Following the strategy in \cite[Section D]{AW2}, one may prove the following result~:

\medskip

\emph{Assume \emph{\textbf{(POT)}} holds and let $s\ge 1$. Then for any bounded $I\subseteq  \sigma_{ac}^{\eps}(\delta)$, we have}
\begin{equation}\label{e:moments}
\sup_{\lambda\in I} \sup_{\eta \in (0,1)} \max_{j=1,\dots,m}\mathbf{E}\left(|\Im \zeta_{j,\eps}^{\lambda+i\eta}|^{-s}\right) \le C_{I,\delta,s} < \infty \, .
\end{equation}

\medskip




This statement essentially means that \textbf{(Green)} holds over $\sigma_{ac}^{\eps}(\delta)$. However, it is not known in general if $\sigma_{ac}^{\eps}(\delta)$ is a large set, so we now explore it further for specific trees.

Following \cite{KLW}, we let
\begin{multline*}
\Sigma = \bigcup \big\{U\subset \R \text{ open}: \forall \,j=1,\dots,m,\ \gamma \mapsto -\zeta_{j,0}^{\gamma} \text{ from } \C^+\to\C^+ \text{ can be extended}\\
\text{to a unique continuous function } \C^+\cup U \to \C^+\big\}\, .
\end{multline*}
By Proposition~\ref{prp:Lang}, if $\T$ satisfies (C1), then $\Sigma \supseteq \sigma(\cA_{\T})\setminus I_{c,C}$.
 
We introduce an additional assumption on $\T$~:

\smallskip

\begin{enumerate}[\textbf{(C2)}]
\item for each $k\in\mathfrak{A}$, there is $k'$ with $M_{k,k'}\ge 1$ such that for any $l\in\mathfrak{A}$~: $M_{k,l}\ge 1$ implies $M_{k',l}\ge 1$.
\end{enumerate}

\smallskip

According to the results of Keller, Lenz and Warzel \cite{KLW}, if conditions \textbf{(C1)} and \textbf{(C2)} hold, then for any compact $I\subset \Sigma$ and $p>1$, we have
$\lim_{\eps \To 0}\sup_{\gamma \in I+i(0,1]} \mathbf{E}(|\zeta_{j,\eps}^{\gamma}-\zeta_{j,0}^{\gamma}|^p) =0$ for any $j$. Combined with the argument in our paper \cite[Proposition 2.2]{AS3}, this implies that $I \subseteq \sigma_{ac}^{\eps}(\delta)$ for some $\delta>0$ if $\eps>0$ is small enough. Summarizing, we have the following.

\medskip

\emph{Assume $\T$ satisfies conditions \emph{\textbf{(C1)}} and \emph{\textbf{(C2)}}. Consider the random Schr\"odinger operator $\cH_{\eps} =\cA+\eps \cW $, with the law of $\cW$ satisfying \emph{\textbf{(POT)}}. Then for any compact $I\subset \sigma(\cA)\setminus I_{c, C}$, the moment bounds \eqref{e:moments} hold on $I$ if $\eps>0$ is small enough.}

\medskip

As an example, consider the tree $\T$ defined by the neighbour matrix $A_{\cN}=\begin{bmatrix} a & b \\ c & d \end{bmatrix}$ with all entries non-zero. Then $(\T,\bullet)$ corresponds to the $5\times 5$ matrix given at the end of \S~\ref{s:Puiseux}, which satisfies \textbf{(C1)} and \textbf{(C2)} if $\min(a,b,c,d)\ge 2$.

\section{Concluding remark}
To complement Theorem \ref{t:thm2}, it would be nice to have some geometric information on the quantity $\langle K\rangle_{\gamma}  = \sum_{x,y\in V_N} K(x,y)\Phi_{\gamma}^N({x},{y})$. How do the weights $\Phi_{\gamma}^N({x},{y})$ depend on $\gamma, x, y$~?
For instance, if $R=0$ and $K(x, x)=a(x)$, we have $\langle a\rangle_{\gamma}  = \sum_{x\in V_N} a(x)\Phi_{\gamma}^N({x},{x})$ where $\Phi_{\gamma}^N({x},{x})$ is a probability measure on $V_N$. For regular graphs,
$\Phi_{\gamma}^N({x},{x})$ is the uniform measure on $V_N$, so it does not depend on $\gamma$ nor on $x$. For irregular graphs, does this probability measure favour vertices with high degree, on on the opposite, low degree~? Does the answer depend on $\gamma$~?

Although a general answer seems far from hand, we can make explicit calculations for the biregular tree $\T_{p,q}$. Note that $\T_{p,q}$ corresponds to the neighbour matrix $A_{\cN}=\begin{bmatrix} 0 & p+1 \\ q+1 & 0 \end{bmatrix}$.
Let $(G_N)$ be a sequence of $(p+1,q+1)$ biregular graphs converging to $\T_{p,q}$. Note that $\T_{p,q}$ is the universal cover of $G_N$ for each $N$.

For simplicity, we only consider $R\le 1$, that is, $K(x,y)$ is supported on the diagonal or on nearest neighbours. For $\lambda$ in the spectrum, one computes that $\frac{\Im G^{\lambda+i0}(\bullet,\bullet)}{\Im G^{\lambda+i0}(\circ,\circ)} = \frac{p+1}{q+1}$, $\frac{\Im G^{\lambda+i0}(\bullet,\bullet)}{\Im G^{\lambda+i0}(\bullet,\circ)}=\frac{p+1}{\lambda}$ and $\frac{\Im G^{\lambda+i0}(\circ,\circ)}{\Im G^{\lambda+i0}(\bullet,\circ)} = \frac{q+1}{\lambda}$.

Recall that $\langle K\rangle_{\lambda+i\eta_0}=\sum_{x,y\in G_N}K(x,y)\frac{\Im \tilde{g}^{\lambda+i\eta_0}_N(x,y)}{\sum_{x\in G_N}\Im \tilde{g}^{\lambda+i\eta_0}_N(x,x)}$, where $\tilde{g}_N^{\gamma}$ is a ``lifted'' Green function. Here, it coincides with the limiting Green function $G^{\gamma}$ of $\T_{p,q}$. In particular, $\sum_{x\in G_N} \Im \tilde{g}_N^{\gamma}(x,x)=\frac{q+1}{p+q+2}N\Im G^{\lambda+i0}(\bullet,\bullet)+\frac{p+1}{p+q+2}N\Im G^{\lambda+i0}(\circ,\circ)$. Thus,
\begin{multline}
\langle K\rangle_{\lambda+i0}
 = \frac{p+q+2}{2N} \\ \cdot \left[\frac{1}{q+1}\sum_{\bullet\in G_N} K(\bullet,\bullet) + \frac{1}{p+1}\sum_{\circ\in G_N} K(\circ,\circ) + \frac{\lambda}{(p+1)(q+1)}\sum_{x\in G_N,y\sim x} K(x,y)\right]\,.
\end{multline}

For $R=0$,  the last term is absent, and we note that $\langle K\rangle_{\lambda+i0}$ does not depend on $\lambda$. From the explicit expression above, we see that the probability measure $\Phi_{\lambda+i0}^N({x},{x})$ gives weight $\frac{p+q+2}{2N(p+1)}$ to each vertex of valency $q+1$, and weight $\frac{p+q+2}{2N(q+1)}$ to each vertex of valency $p+1$. The total measure of vertices of valency $p+1$ goes to $1/2$ as $N\To +\infty$, and so does the total mass of vertices of valency $q+1$.

We see the influence of $\lambda$ on the term corresponding to $R=1$.

\medskip

{\bf{Acknowledgements~:}} This material is based upon work supported by the Agence Nationale de la Recherche under grant No.ANR-13-BS01-0007-01, by the Labex IRMIA and the Institute of Advance Study of Universit\'e de Strasbourg, and by Institut Universitaire de France.

\providecommand{\bysame}{\leavevmode\hbox to3em{\hrulefill}\thinspace}
\providecommand{\MR}{\relax\ifhmode\unskip\space\fi MR }
\providecommand{\MRhref}[2]{%
  \href{http://www.ams.org/mathscinet-getitem?mr=#1}{#2}
}
\providecommand{\href}[2]{#2}

\end{document}